\documentclass[11pt,a4]{article}
\usepackage{latexsym}
\usepackage{multirow}
\usepackage{amssymb,amsmath}
\usepackage{graphicx}
\input epsf
\makeatletter
\def\section{\@startsection {section}{1}{\z@}{-3.5ex plus -1ex minus
-.2ex}{2.3ex plus .2ex}{\large\sc}}
\def\subsection{\@startsection{subsection}{2}{\z@}{-3.25ex plus -1ex minus
 -.2ex}{1.5ex plus .2ex}{\normalsize\sc}}
\makeatother
\makeatletter \@addtoreset{equation}{section}

\makeatother
\newcommand{\nc}{\newcommand}

\nc{\bea}{\begin{eqnarray}} \nc{\eea}{\end{eqnarray}}
\nc{\be}{\bea} \nc{\ee}{\eea}
\nc{\tr}{\mathop{\mbox{tr}}\nolimits}
\nc{\ad}{\mathop{\mbox{ad}}\nolimits}
\nc{\Tr}{\mathop{\mbox{Tr}}\nolimits}
\nc{\Det}{\mathop{\mbox{Det}}\nolimits}
\nc{\rk}{\mathop{\mbox{rk}}\nolimits} \nc{\ra}{\rightarrow}
\nc{\Ra}{\Rightarrow} \nc{\LRa}{\Leftrightarrow} \nc{\ot}{\otimes}
\nc{\non}{\nonumber\\} \nc{\ZZ}{\mathbb{Z}} \nc{\RR}{\mathbb{R}}

\newtheorem{theorem}{Theorem}
\newtheorem{proposition}{Proposition}
\newtheorem{lemma}{Lemma}
\newtheorem{corollary}{Corollary}

\newtheorem{remark}{Remark}
\newtheorem{question}{Question}
\def\1#1{^{(#1)}}
\def\p{\prime}

\def\ra{\rangle}
 
\begin{document}
\title{Summatory Multiplicative Arithmetic Functions:
\\Scaling and Renormalization}  
\author{Leonid G. Fel\\
Department of Civil Engineering, Technion, Haifa 32000, Israel}
\vspace{-.2cm}
\date{}
\maketitle
\def\be{\begin{equation}}
\def\ee{\end{equation}}
\def\bea{\begin{eqnarray}}
\def\eea{\end{eqnarray}}
\def\p{\prime}
\begin{abstract}
We consider a wide class of summatory functions $F\left\{f;N,p^m\right\}=\sum_{
k\leq N}f\left(p^m k\right)$, $m\in{\mathbb Z}_+\cup\{0\}$, associated with the 
multiplicative arithmetic functions $f$ of a scaled variable $k\in{\mathbb Z}_+
$, where $p$ is a prime number. Assuming an asymptotic behavior of summatory
function, $F\{f;N,1\}\stackrel{N\to\infty}{=}G_1(N)\left[1+{\cal O}\left(G_2(N)
\right)\right]$, where $G_1(N)=N^{a_1}\left(\log N\right)^{b_1}$, $G_2(N)=N^{-a
_2}\left(\log N\right)^{-b_2}$ and $a_1,a_2\geq 0$, $-\infty<b_1,b_2<\infty$, 
we calculate a renormalization function defined as a ratio, $R\left(f;N,p^m
\right)=F\left\{f;N,p^m\right\}/F\{f;N,1\}$, and find its asymptotics $R_{
\infty}\left(f;p^m\right)$ when $N\to\infty$. We prove that the renormalization 
function is multiplicative, i.e., $R_{\infty}\left(f;\prod_{i=1}^np_i^{m_i}
\right)=\prod_{i=1}^nR_{\infty}\left(f;p_i^{m_i}\right)$ with $n$ distinct 
primes $p_i$. We extend these results on the others summatory functions $\sum_{
k\leq N}f(p^mk^l)$, $m,l,k\in{\mathbb Z}_+$ and $\sum_{k\leq N}\prod_{i=1}^nf_i
\left(kp^{m_i}\right)$, $f_i\neq f_j$, $m_i\neq m_j$. We apply the derived 
formulas to a large number of basic summatory functions including the Euler 
$\phi(k)$ and Dedekind $\psi(k)$ totient functions, divisor $\sigma_n(k)$ and 
prime divisor $\beta(k)$ functions, the Ramanujan sum $C_q(n)$ and Ramanujan 
$\tau$ Dirichlet series, and others.
\end{abstract}
\tableofcontents
\section{Summatory Multiplicative Functions with Scaled Variable}\label{r1}
Among summatory arithmetic functions $\sum_{k\leq N}f(k),\;k\in{\mathbb Z}_+$, 
of various $f(k)$ the most utilized are the basic multiplicative functions 
$f(k)$ and their algebraic combinations. Study of different summatory functions 
and their asymptotics has a long history \cite{apo95}, \cite{dic50} and 
\cite{har79}. Their list includes totient functions: the Euler $\phi(k)$, 
Dedekind $\psi(k)$ and Jordan $J_n(k)$, and non totient functions: the M\"obius 
$\mu(k)$ and the $n$-order M\"obius $\mu_n(k)$, Liouville $\lambda(k)$, Piltz 
$d_n(k)$, divisor $\sigma_n(k)$, prime divisor $\beta(k)$, non isomorphic 
Abelian group enumeration function $\alpha(k)$, exponentiation of additive 
functions $\omega(k)$ and $\Omega(k)$ which give the numbers of distinct prime 
dividing $k$ and total prime factors of $k$ counted with multiplicities, 
respectively. The whole family of multiplicative arithmetic functions is much 
wider, e.g., the number $q_n(k)$ of representations of $k$ by sum of two 
integral $n$th powers \cite{kr88}, the number $r_n(k)$ of representations of $k$
by sum of $n$ integer squares \cite{har79}, the Legendre and Zsigmondy totient 
functions \cite{dic50} and the Nagell totient function \cite{hor68}, the non 
isomorphic solvable \cite{gn70} and nilpotent \cite{lub03} finite group 
enumeration functions, the Gauss \cite{har79}, Ramanujan \cite{har79} and 
Kloosterman \cite{gs83} sums, the Ramanujan $\tau(k)$ function \cite{leh43} and 
others.

In this article we study a family of summatory multiplicative arithmetic 
functions a with scaled summation variable, $F\left\{f;N,p^m\right\}=\sum_{k
\leq N}f(p^mk),\;k,m\in{\mathbb Z}_+$, where $p$ is a prime number. For this 
reason we use hereafter the notation $F\{f;N,1\}$ for unscaled summatory 
function. A description of asymptotics of $F\{f;N,1\}$, $N\to\infty$, assumes 
that we know two characteristics, its leading and error terms, $G_1(N)$ and 
$G_2(N)$, i.e., $F\{f;N,1\}\stackrel{N\to\infty}{=}G_1(N)\left[1+{\cal O}\left(
G_2(N)\right)\right]$. 

In section \ref{r11} we introduce universality classes ${\mathbb B}\{G_1(N);G_2(
N)\}$ of arithmetic functions $f(k)$ such that different functions possess the 
same $G_1(N)$ and $G_2(N)$. By inspection of a vast number of multiplicative 
functions $f(k)$ we focus on their wide class, $G_1(N)=N^{a_1}\left(\log N
\right)^{b_1}$, $G_2(N)=N^{-a_2}\left(\log N\right)^{-b_2}$, where $a_1,a_2\geq 
0$ and $-\infty<b_1,b_2<\infty$. In section \ref{r12} we derive a functional 
equation defined at different scales,
\bea
F\left\{f;N,p^m\right\}=\sum_{r=0}^{\left\lfloor\log_p N\right\rfloor}L_r\left(
f;p^m\right)F\left\{f;\left\lfloor\frac{N}{p^r}\right\rfloor,1\right\}\;,
\label{i0}
\eea
where the characteristic functions $L_r\left(f;p^m\right)$ are satisfied the 
recursive equations,
\bea
L_r\left(f;p^m\right)=f\left(p^{m+r}\right)-\sum_{j=0}^{r-1}L_j\left(f;p^m
\right)f\left(p^{r-j}\right)\;,\quad L_0\left(f;p^m\right)=f\left(p^m\right)\;,
\quad f(1)=1\;,\label{i1}
\eea
and $\left\lfloor u\right\rfloor$ denotes the largest integer not exceeding $u$.
The functions $L_r\left(f;p^m\right)$ are calculated in (\ref{b6}) and their 
behavior in $r$ is crucial for convergence of numerical series. This is a 
subject of special discussion in the next section.

In section \ref{r2} we define the renormalization function $R\left(f;N,p^m
\right)$ and its asymptotics,
\bea
R\left(f;N,p^m\right)=\frac{F\left\{f;N,p^m\right\}}{F\left\{f;N,1\right\}}\;,
\quad R_{\infty}\left(f,p^m\right)=\lim_{N\to\infty}R\left(f;N,p^m\right)\;.
\label{i2}
\eea
The aim of this paper is to study the asymptotic renormalization function $R_{
\infty}\left(f,p^m\right)$ in various aspects: (a) its existence as the 
convergent numerical series, (b) its multiplicativity property without 
specifying the 
function $f(k)$, (c) formulas for $R_{\infty}\left(f(k)\cdot k^{-s},p^m\right)$ 
for corresponding Dirichlet series, (d) formulas for $R_{\infty}\left(f(k^n),
p^m\right)$ for different arithmetic functions $f$, (e) formulas $R_{\infty}
\left(f,p^m\right)$ for basic arithmetic function $f$. For short, we'll often 
skip the word 'asymptotic' and refer to $R_{\infty}\left(f,p^m\right)$ as 
renormalization function if this would not mislead the readers.

Imposing the constraints on $L_r\left(f;p^m\right)$, in section \ref{r21} we
prove two Lemmas on convergence of numerical series and calculate asymptotics 
of renormalization functions. In section \ref{r22} we show that by these 
constraints the error term $G_2(N)$ does not contribute to $R_{\infty}\left(f,
p^m\right)$. In section \ref{r23} we give a rational representation for $R_{
\infty}\left(f,p^m\right)$ which is much easier to implement in analytic
calculations.

In section \ref{r3} we prove that the renormalization function has 
multiplicative property in the following sense, $R_{\infty}\left(f;\prod_{i=1}^n
p_i^{m_i}\right)=\prod_{i=1}^nR_{\infty}\left(f;p_i^{m_i}\right)$ with $n$ 
distinct primes $p_i\geq 2$. Making use of renormalization approach we also 
calculate the summatory functions $\sum_{k_1,k_2\leq N}f(k_1k_2)$.

In section \ref{r4} we extend the renormalization approach on summatory $\sum_{
k\leq N}\prod_{i=1}^nf_i\left(kp^{m_i}\right)$, where $f_i\neq f_j$, $m_i\neq 
m_j$, and the corresponding Dirichlet series $\sum_{k=1}^{\infty}f(kp^m)k^{-s}$,
$k,n\in{\mathbb Z}_+$. We also study the renormalization of summatory function 
$\sum_{k\leq N}f\left(k^np^m\right)$.

In section \ref{r5} we apply formulas, derived in section \ref{r2}, to calculate
$R_{\infty}\{f;p^m\}$ for basic multiplicative arithmetic functions $f$ and 
their combinations. Almost all summatory functions are treated by Theorem 
\ref{thm2} based on a simple calculation of $f\left(p^r\right)$ and avoiding a 
cumbersome calculation of characteristic functions $L_r\left(f;p^m\right)$. The 
renormalization functions are given by algebraic and non-algebraic expressions 
as well, e.g., see $R_{\infty}\left(\sigma_0^n;p^m\right)$ in (\ref{r521d}) and 
$R_{\infty}\left(\sigma_1/\sigma_0;p\right)$ in (\ref{r521g}), respectively. In 
this conjunction, the Ramanujan $\tau$ function is of particular interest: in 
contrast to many other functions $f(k)$ its value at $k=p^r$ is given by a heavy
formula (\ref{r531}), while the characteristic functions $L_r\left(\tau;p^m
\right)$ and $L_r\left(\tau^2;p^m\right)$ have been calculated in a simple form 
suitable for explicit calculation of $R_{\infty}\left(\tau\cdot k^{-s},p^m
\right)$ and $R_{\infty}\left(\tau^2,p^m\right)$. We have found a new identity 
(\ref{r538}) for the Ramanujan $\tau$ function.

In section \ref{r6} we give a numerical verification to renormalization approach
developed in this article by numerical calculations and show its validity with 
high precision.
\subsection{Asymptotic Growth of Summatory Functions}\label{r11}
Consider the summatory multiplicative arithmetic function $F\{f;N,1\}$ and 
represent its asymptotics in $N$ by using one constant ${\cal F}$ and two 
positive definite functions $G_1(N)$ and $G_2(N)$,
\bea
F\{f;N,1\}\stackrel{N\to\infty}{=}{\cal F}\;G_1(N)\left[1+{\cal O}\left(G_2(N)
\right)\right]\;,\quad\lim_{N\to\infty}G_2(N)=0\;,\quad 
G_1(N),\;G_2(N)>0\;,\hspace{.5cm}\label{i4}
\eea
where $"{\cal O}"$ stands for the "big--{\cal O}" Landau symbol. 

Asymptotic growth of $F\{f;N,1\}$ is determined by its leading term and is given
by nondecreasing function $G_1(N)$, i.e., either increasing or unity, while the 
decreasing function $G_2(N)$ stands for the error term. The constant ${\cal F}$ 
is introduced to distinguish summatories with similar functions $G_1(N)$ and
$G_2(N)$, e.g., see $F\{1/\phi;N,1\}$, $F\{1/\psi;N,1\}$ and $F\{1/\sigma_1;N,
1\}$ in Table 1.

Different arithmetic functions $f(k)$ may possess the same $G_1(N)$ and $G_2(
N)$, therefore the whole set of $f(k)$ can be decomposed into different 
universality classes ${\mathbb B}\{G_1(N);G_2(N)\}$ as follows,
\bea
f(k)\in{\mathbb B}\{G_1(N);G_2(N)\}\;,\quad\mbox{where}\label{i5}
\eea
\bea
{\mathbb B}\{G_1(N);G_2(N)\}\!=\!\left\{f(k)\left|\;\frac{F\{f;N,1\}}{{\cal F}\;
G_1(N)}\right.\!\stackrel{N\to\infty}{=}\!1,\;\left|\frac{F\{f;N,1\}}
{{\cal F}\;G_1(N)}-1\right|\!\stackrel{N\to\infty}{<}\!{\cal C}G_2(N);\;0<{\cal 
C},{\cal F}<\infty\right\}\nonumber
\eea
Below we give examples of various multiplicative functions $f(k)$ which belong 
to the different universality classes,
\bea
\frac1{\phi(k)},\;\frac1{\psi(k)},\;\frac{\phi(k)}{k^2},\;\frac1{\sigma_1(k)}
\in{\mathbb B}\left\{\log N;\frac1{\log N}\right\},\quad \frac{\phi(k)}{k},\;
\frac{\sigma_1(k)}{k},\;\sigma_{-1}(k)\in{\mathbb B}\left\{N;\frac{\log N}{N}
\right\},\nonumber\\
\sigma_0(k),\;2^{\omega(k)},\;\mu^2(k)2^{\omega(k)}\in{\mathbb B}\left\{N\log 
N;\frac1{\log N}\right\},\quad \frac{\mu^2(k)}{k},\;\alpha(k),\;\beta(k)\in
{\mathbb B}\left\{N;\frac1{\sqrt{N}}\right\}.\nonumber
\eea
By inspection of a vast number of multiplicative functions $f(k)$ with known 
asymptotics $G_1(N)$ and $G_2(N)$ (see Tables 1 and 2) in this paper we focus on
their most wide class,
\bea
G_1(N)=N^{a_1}\left(\log N\right)^{b_1},&&\left\{\begin{array}{lr}a_1>0,&-
\infty<b_1<\infty\;,
\\a_1=0,&0\leq b_1<\infty\;,\end{array}\right.\label{i7}\\
G_2(N)=\frac1{N^{a_2}\left(\log N\right)^{b_2}},&&\left\{\begin{array}{lr}a_2>0,
&-\infty<b_2<\infty\;,\\a_2=0,&0<b_2<\infty\;.\end{array}\right.\nonumber
\eea
\begin{remark}\label{rem1}Tables 1 and 2 do not present any example of 
multiplicative arithmetic functions $f(k)$ of the special universality classes 
${\mathbb B}\left\{N^{a_1}\left(\log N\right)^{b_1};\;N^{-a_2}\left(\log N
\right)^{-b_2}\right\}$ such that
\bea
0<a_1<a_2\;,\;-\infty<b_1,b_2<\infty\hspace{.5cm}\mbox{and}\hspace{.5cm}
a_1=0<a_2\;,\;0\leq b_1<b_2<\infty\;.\label{i6}
\eea
Despite an extensive search in the available literature we have not found 
such functions there.
\end{remark}
In Tables 1 and 2 we present a long list of multiplicative functions which
belong to one of the universality classes ${\mathbb B}\left\{G_1(N);G_2(N)
\right\}$ satisfying (\ref{i7}). We use the standard notations for the functions
mentioned in section \ref{r1}. Here $\zeta(s)$ stands for the Riemann zeta
function, the explicit expressions for $A_n$, $B_s$, $D_n$, $E_n$, $K_n$ and
$I_n$ and values for $C_i$ are given in the corresponding references. The values
of $C_2$, $C_8$ and $C_{10}$ were calculated by author and marked by ($\star$).
The error terms for summatories of $|\tau(k)|$, $\tau^2(k)$, $\tau^4(k)$ and
$\tau^2(k)/\sqrt{k^{25}}$ are unknown to date.

\begin{center}
{\bf Table$\;$1}. Summatory multiplicative arithmetic functions and their
asymptotics.\\
\vspace{.1cm}
\begin{tabular}{|c||c|c|c|c|c|c|c|c|c|c|c|c||} \hline
$f(k)$ & ${\cal F}$ & $G_1(N)$ & $G_2(N)$  & Ref \\ \hline\hline
$\mu^2(k)$ & $\zeta^{-1}(2)$ & $N$ & $N^{-1/2}$ 
& \cite{br02}\\ \hline
$\mu^2(k)/\phi(k)$ & 1 & $\log N$ & $\log^{-1}N$ 
& \cite{war27}\\ \hline
$\mu_n(k)$ & $A_n$, $n\geq 2$ & $N$ & $N^{\frac1{n}-1}\log N$ 
& \cite{san96}, p.193\\ \hline
$1/\phi(k)$ & $\zeta(2)\zeta(3)/\zeta(6)$ & $\log N$ & $\log^{-1}N$
& \cite{lan00}\\ \hline
$k/\phi(k)$ & $\zeta(2)\zeta(3)/\zeta(6)$ & $N$ & $N^{-1}\log N$ 
& \cite{cit85}\\ \hline
$\psi(k)$ & $1/2\;\zeta(2)/\zeta(4)$ & $N^2$ &$N^{-1}\log N$
& \cite{apo95}, p.72\\ \hline
$1/\psi(k)$ & $C_1\simeq 0.37396$ & $\log N$ & $\log^{-1}N$ 
& \cite{sita79} \\ \hline
$\sigma_{-b}(k)$ & $\zeta(b+1)$ & $N$, $b>0$, $b\neq 1$ & $N^{-min\{1,b\}}$ 
& \cite{apo95}, p.61 \\ \hline
$\sigma_{-1}(k)$ & $\zeta(2)$ & $N$ & $N^{-1}\log N$ 
& \cite{apo95}, p.61 \\ \hline
$\sigma_a(k)$ & $\zeta(a+1)/(a+1)$&$N^{a+1}$, $a>0$, $a\neq 1$&$N^{-min\{1,a\}}$
& \cite{apo95}, p.60\\ \hline
$\sigma_1(k)$ & $\zeta(2)/2$ & $N^2$ & $N^{-1}\log^{2/3} N$ 
& \cite{wal63} \\ \hline
$\sigma_a^2(k)$ & $\frac{\zeta^2(1+|a|)\;\zeta(1+2|a|)}{\zeta(2+2|a|)}$ 
& $N^{1+a+|a|}$, $0<|a|<1$ & $N^{-|a|}\log N$ 
& \cite{pet98} \\ \hline
$\sigma_1^2(k)$ & $5/6\;\zeta(3)$ & $N^3$ & $N^{-1}\log^{5/3} N$ 
& \cite{sm70} \\\hline
$\sigma_0^n(k)$ & $D_n$,$D_1\!=\!1$,$D_2\!=\!\pi^{-2}$ & $N\cdot\log^{2^n-1}N$
& $\log^{-1}N$ 
& \cite{ram15}, \cite{wil23} \\ \hline
$1/\sigma_0(k)$ & $D_{-1}\simeq 0.5469$ & $N\cdot\log^{-1/2}N$ & $\log^{-1}N$ 
& \cite{ram15}, \cite{wil23} \\\hline
$\sigma_1(k)/\phi(k)$ & $C_2\simeq 3.6174\quad\star$ & $N$ & $N^{-1}\log^2N$
& \cite{pet98}\\\hline
$\sigma_1(k)/\sigma_0(k)$ &$C_3\simeq 0.3569$ &$N^2\log^{-1/2} N$& $\log^{-1}N$
& \cite{pom81} \\ \hline
$1/\sigma_1(k)$ & $C_4\simeq 0.6728$ & $\log N$ & $\log^{-1}N$  
& \cite{sita79} \\\hline
$d_n(k)$ & $1/\Gamma(n)$ & $N\cdot\log^{n-1}N$ & $\log^{-1}N$ 
& \cite{har15}, \cite{sel54}\\\hline
$d_n^2(k)$ & $E_n$, $E_2=D_2$ & $N\cdot\log^{n^2-1}N$ & $\log^{-1}N$
& \cite{kal74} \\\hline
$1/d_n(k)$ & $K_n$, $K_2=D_{-1}$ & $N\cdot\log^{1/n-1}N$ &$\log^{-1}N$
& \cite{ivi77} \\\hline
$\beta(k)$ & $\zeta(2)\zeta(3)/\zeta(6)$ &$N$& $N^{-1/2}$
& \cite{ken47}\\ \hline
$\alpha(k)$ & $\prod_{l=2}^{\infty}\zeta(l)\simeq 2.29486$ & $N$ & $N^{-1/2}$ 
& \cite{es35}\\ \hline
$1/\alpha(k)$ & $C_5\simeq 0.75204$ & $N$ &$N^{-1/2}\log^{-1/2} N$
& \cite{kon80}, p.16 \\ \hline
$2^{\omega(k)}$& $\zeta^{-1}(2)$ &$N\cdot\log N$&$\log^{-1}N$
& \cite{sel54}\\ \hline
$3^{\omega(k)}$& $C_6\simeq 0.14338$ & $N\cdot\log^2N$&$\log^{-1}N$ 
& \cite{ten95}, p.53\\ \hline
$2^{\Omega(k)}$ & $C_7\simeq 0.27317$ &$N\cdot\log^2N$&$\log^{-1}N$
& \cite{gro56} \\\hline
$q_n(k)$ & $2\Gamma^2\!\left(n^{-1}\right)\!/\!\!\left(n\Gamma\!\left(2n^{-1} 
\right)\!\right)$ & $N^{2/n}$, $n\geq 3$ & $N^{-1/n(1+1/n)}$ 
& \cite{kr88}. p.143\\\hline
$r_2(k)\!=\!q_2(k)$ & $\pi$  & $N$ & $N^{-1/2}$ & \cite{har79} \\\hline
$|\tau(k)|$ & $C_8\simeq 0.0996\quad\star$ & $N^{13/2}\log^{-1+8/(3\pi)}N$ & ? 
& \cite{ell80}, \cite{ran83} \\ \hline\
$\tau^2(k)$ & $C_9\simeq 0.032007\;$ & $N^{12}$  & ? 
& \cite{leh43}\\ \hline
$\tau^4(k)$ & $C_{10}\simeq 0.0026\quad\star$ & $N^{23}\log N$ & ? 
& \cite{ran83} \\ \hline
\end{tabular}
\label{ta1}
\end{center}
\begin{center}  
{\bf Table$\;$2}. Other summatory functions and the Dirichlet series of 
multiplicative functions.\\
\vspace{.5cm}
\begin{tabular}{|c||c|c|c|c|c|c|c|c|c|c|c|c||} \hline
$f(k)$ & ${\cal F}$ & $G_1(N)$ & $G_2(N)$  & Ref \\ \hline\hline
$\mu(k)/k^s$ & $\zeta^{-1}(s)$, $s>1$ & $N^0$ & $N^{-(s-1)}$ 
& \cite{apo95}, p.231 \\ \hline
$\mu^2(k)/k^s$ & $\zeta(s)/\zeta(2s)$, $s>1$ & $N^0$ & $N^{-(s-1)}$
& \cite{apo95}, p.241 \\ \hline
$\mu^2(k)/k$ & $\zeta^{-1}(2)$ & $\log N$ & $\log^{-1}N$
& \cite{san06}, p.195\\ \hline
$\lambda(k)/k^s$ & $\zeta(2s)/\zeta(s)$, $s>1$ &$N^0$ & $N^{-(s-1)}$ 
& \cite{apo95}, p.231 \\ \hline
$\phi(k)/k^s$ & $\zeta^{-1}(2)/(2-s)$, $0<s\leq 1$ & $N^{2-s}$ & $N^{-1}\log N$
& \cite{apo95}, p.71 \\ \hline
$\phi(k)/k^s$ & $\zeta^{-1}(2)/(2-s)$, $1<s<2$ &$N^{2-s}$ &$N^{-(2-s)}$
& \cite{apo95}, p.71 \\ \hline
$\phi(k)/k^2$ & $\zeta^{-1}(2)$ & $\log N$&$\log^{-1}N$
&\cite{apo95}, p.71 \\ \hline   
$\phi(k)/k^s$ & $\zeta(s-1)/\zeta(s)$, $s>2$ &$N^0$ &$N^{-(s-2)}$
& \cite{apo95}, p.231 \\ \hline
$\left(\phi(k)/k\right)^s$ & $B_s$, $B_1\!=\!\zeta^{-1}(2)$ & $N$, $s>0$ &
$N^{-1}\log^s N$
& \cite{chl27}\\ \hline
$\sigma_0(k)/k^s$ & $(1-s)^{-1}$, $0<s<1$ & $N^{1-s}\log N$ & $\log^{-1}N$ 
& \cite{apo95}, p.70 \\\hline
$\sigma_0(k)/k$ & 1/2 & $\log^2 N$ & $\log^{-1}N$ 
& \cite{apo95}, p.70 \\\hline 
$\sigma_0(k)/k^s$ & $\zeta^2(s)$, $s>1$ & $N^0$ & $N^{-(s-1)}\log N$ 
& \cite{apo95}, p.231 \\\hline
$(\sigma_1(k)/k)^s$ & $I_s$, $I_1=\zeta(2)$, $I_2=\frac{5}{2}\zeta(3)$
& $N$, $s>0$ & $N^{-1}\log^s N$
& \cite{sm70}, \cite{pet98}\\ \hline
$r_2(k)/k$ & $\pi$ & $\log N$ & $\log^{-1}N$ & \cite{har79} \\\hline
$\tau^2(k)/k^{25/2}$ & $C_{11}\simeq 1.58824$ & $N^0$  & ?
& \cite{do04} \\ \hline
\end{tabular}
\label{ta2}
\end{center}
A relationship between the multiplicative properties of arithmetic functions and
asymptotics of their summatory functions is not straightforward. In other words,
a correspondence:
\bea
f(k)\;\;\mbox{is a multiplicative function}\;\;\;\longleftrightarrow\;\;\;f(k)
\in{\mathbb B}\left\{N^{a_1}\left(\log N\right)^{b_1};\;N^{-a_2}\left(\log N
\right)^{-b_2}\right\}\label{i8}
\eea
is neither bijective nor injective. Indeed, the direction `$\longleftarrow$` is 
not holding since there exists a non multiplicative function $f(k)=\log\phi(k)/
\log\sigma_1(k)$ which has the summatory function $F\left\{f;N,1\right\}
\stackrel{N\to\infty}=N+{\cal O}\left(N\log^{-1}N\right)$ \cite{kon80}.

Regarding another direction `$\longrightarrow$` there exist multiplicative 
arithmetic functions with summatory growth that differs from $N^{a_1}\left(\log 
N\right)^{b_1}$ and come by enumeration of finite groups. Let $\chi(k)$ be a 
number of nilpotent groups of order $k$ which is multiplicative, because each 
finite nilpotent group is a direct product of its Sylow subgroups \cite{lub03}. 
When $k=p^r$, for a prime $p$, it is known \cite{hig60} that $\chi(k)\simeq p^{
(2/27+o(1))\;r^3}$, i.e., of the order $k^{(\log_pk)^2}$. Consider its summatory
$F[\chi;N,1]$ which is the number of nilpotent groups of order at most $N$. For 
some $r$ we have $2^r\le N<2^{r+1}$. Then for some $0<s<1$ we have,
\bea
F[\chi;N,1]\ge \chi(2^r)=2^{(2/27+o(1))\;r^3}\ge (2^{r+1})^{2s/27\;(r+1)^2}\;
\;\rightarrow\;\; F[\chi;N,1]\ge N^{2s/27\;(\log N)^2},\label{i9}
\eea
and therefore, $\chi(k)\not\in{\mathbb B}\left\{N^{a_1}\left(\log N\right)^{b_1}
;\;N^{-a_2}\left(\log N\right)^{-b_2}\right\}$. Note that the number $\chi(k,c,
d)$ of nilpotent groups of order $k$, of nilpotency class at most $n$, generated
by at most $m$ elements does belong to the universality class defined in 
(\ref{i7}) where $a_i$ and $b_i$ depend on $n$ and $m$ \cite{sau03}.

Other examples of multiplicative function with summatory growth $N^{a_1}\left(
\log N\right)^{b_1}$ and the error term, which is different than $N^{-a_2}\left(
\log N\right)^{-b_2}$, were given in \cite{wal63},
\bea
\quad\phi(k)\in{\mathbb B}\left(N^2;\frac1{N}\left(\log N\right)^{2/3}\left(
\log\log N\right)^{4/3}\right),\quad\frac{\phi(k)}{k}\in{\mathbb B}\left(N;
\frac1{N}\left(\log N\right)^{2/3}\left(\log\log N\right)^{4/3}\right).\nonumber
\eea
They can also be encompassed within the universality classes by extending the 
latter on much wider family of asymptotics, e.g., ${\mathbb B}\left\{N^{a_1}
\left(\log N\right)^{b_1}\left(\log\log N\right)^{c_1};\;N^{-a_2}\left(\log N
\right)^{-b_2}\left(\log\log N\right)^{-c_2}\right\}$. Keeping in mind such 
option, we continue to study the summatory multiplicative functions with 
universality classes of asymptotics given in (\ref{i7}).
\subsection{Scaling Equation for Summatory Functions}\label{r12}
Represent $F\left\{f;N,p^m\right\}$ as a sum of two summatory functions
\bea
F\left\{f;N,p^m\right\}=\sum_{k=1\atop p\;\nmid\;k}^Nf\left(p^mk\right)+
\sum_{k=p\atop p\;\mid\;k}^Nf\left(p^mk\right)\;.\label{b1}
\eea
Making use of multiplicativity, $f(k_1k_2)=f(k_1)f(k_2)$ if $\gcd(k_1k_2)=1$, 
$f(1)=1$, we get
\bea
F\left\{f;N,p^m\right\}\!\!\!&=&\!\!\!f\left(p^m\right)\sum_{k=1\atop p\;\nmid\;
k}^Nf(k)+\sum_{k=p\atop p\;\mid\;k}^Nf\left(p^mk\right)=f\left(p^m\right)\left(
\sum_{k=1}^Nf(k)-\sum_{k=p\atop p\;\mid\;k}^Nf(k)\right)+\sum_{k=p\atop p\;\mid
\;k}^Nf\left(p^mk\right)\nonumber\\
&=&\!\!\!f\left(p^m\right)F\{f;N,1\}-f\left(p^m\right)\sum_{l=1}^{N_1}f(pl)+
\sum_{l=1}^{N_1}f\left(p^{m+1}l\right)\nonumber\\
&=&\!\!\!f\left(p^m\right)F\{f;N,1\}-f\left(p^m\right)F\{f;N_1,p\}+
F\left\{f;N_1,p^{m+1}\right\}.\label{b2}
\eea
where $N_r=\left\lfloor N/p^r\right\rfloor$. The recursion (\ref{b2}) holds for 
any $N_r$, i.e.,
\bea
F\left\{f;N_r,p^m\right\}=f\left(p^m\right)F\{f;N_r,1\}-f\left(p^m\right)
F\{f;N_{r+1},p\}+F\left\{f;N_{r+1},p^{m+1}\right\}\;.\label{b3}
\eea
Substituting (\ref{b3}) into (\ref{b2}) we obtain
\bea
F\left\{f;N,p^m\right\}\!\!&=&\!\!f\left(p^m\right)F\{f;N,1\}-f\left(p^m\right)
\left[f(p)F\{f;N_1,1\}-f(p)F\{f;N_2,p\}+F\left\{f;N_2,p^2\right\}\right]
\nonumber\\
&+&\!\!f\left(p^{m+1}\right)F\{f;N_1,1\}-f\left(p^{m+1}\right)F\{f;N_2,p\}+
F\left\{f;N_2,p^{m+2}\right\}\nonumber\\
&=&\!\!f\left(p^m\right)F\{f;N,1\}+\left[f\left(p^{m+1}\right)-f\left(p^m
\right)f(p)\right]F\{f;N_1,1\}\nonumber\\
&+&\!\!\left[f\left(p^m\right)f(p)-f\left(p^{m+1}\right)\right]F\{f;N_2,p\}-
f\left(p^m\right)F\left\{f;N_2,p^2\right\}+F\left\{f;N_2,p^{m+2}\right\}.
\nonumber
\eea
Continuing this procedure recursively we get finally,
\bea
&&\hspace{-.5cm}
F\left\{f;N,p^m\right\}=\sum_{r=0}^{\left\lfloor\log_p N\right\rfloor}
L_r\left(f;p^m\right)F\{f;N_r,1\}\;,\quad\mbox{where}\label{b4}\\
&&\hspace{-.5cm}
L_r\left(f;p^m\right)=f\left(p^{m+r}\right)-\sum_{j=0}^{r-1}L_j\left(f;p^m
\right)f\left(p^{r-j}\right)\;,\quad 
L_0\left(f;p^m\right)=f\left(p^m\right).\label{b5}
\eea
The straightforward calculations of $L_r\left(f;p^m\right)$ give
\bea
L_1\left(f;p^m\right)\!\!\!&=&\!\!\!
f\left(p^{m+1}\right)-f\left(p^m\right)f(p),\label{b6}\\
L_2\left(f;p^m\right)\!\!\!&=&\!\!\!f\left(p^{m+2}\right)-f\left(p^m\right)
f\left(p^2\right)-f\left(p^{m+1}\right)f(p)+f\left(p^m\right)f^2(p),\nonumber\\
L_3\left(f;p^m\right)\!\!\!&=&\!\!\!f\left(p^{m+3}\right)-f\left(p^m\right)
f\left(p^3\right)-f\left(p^{m+1}\right)f\left(p^2\right)+2f\left(p^m\right)
f\left(p^2\right)f(p)-\nonumber\\
&&\!\!\!f\left(p^{m+2}\right)f(p)+f\left(p^{m+1}\right)f^2(p)-
f\left(p^m\right)f^3(p),\quad\mbox{etc}\;,\nonumber
\eea
such that for $p=1$ or $m=0$ we have, $L_0\left(f;1\right)=1$ and $L_r\left(f;
1\right)=0$, $r\geq 1$. 

By (\ref{b5}) or (\ref{b6}) the general formulas for $L_r\left(f;p^m\right)$ 
can be calculated by induction for simple arithmetic functions $f(k)$ and 
$1/f(k)$ such that $f(p^m)=A_f\;p^{m-1}$, $m\geq 1$,
\bea
L_r\left(f;p^m\right)=A_f(p-A_f)^rp^{m-1},\quad L_r\left(\frac1{f};p^m\right)=
\frac1{A_f}\left(\frac1{p}-\frac1{A_f}\right)^r\frac1{p^{m-1}},\quad r\geq 
1,\quad\label{b6a}
\eea
and $A_f$ denotes a real constant, e.g., $A_{\phi}=p-1$, $A_{\psi}=p+1$ that 
gives
\bea
\left\{\begin{array}{l}L_r\left(\phi;p^m\right)=p^{m-1}(p-1),\\
L_r\left(\psi;p^m\right)=(-1)^rp^{m-1}(p+1),\end{array}\right.\;
L_r\left(\frac1{\phi};p^m\right)=\frac{(-1)^rp^{1-r-m}}{(p-1)^{r+1}}\;,\;\;
L_r\left(\frac1{\psi};p^m\right)=\frac{p^{1-r-m}}{(p+1)^{r+1}}.\nonumber
\eea
Another example of arithmetic functions leading to $L_r\left(f;p^m\right)=0$, 
$r\geq 1$, is that when $A_f=p$ in (\ref{b6a}) or $f(p^m)=c^m$, e.g., $2^{
\Omega(p^m)}=2^m$ and $\lambda(p^m)=(-1)^m$.

In general case of $f(k)$ the formula of $L_r\left(f;p^m\right)$ with arbitrary 
$r\geq 0$ can be hardly recognized by its partial expressions, e.g., for $f(k)=
1/\sigma_0(k)$,
\bea
L_0=\frac1{m+1},\quad L_1=\frac{m\;m!}{2(m+2)!},\quad L_2=\frac{(5m+7)\;m\;m!}
{12(m+3)!},\quad L_3=\frac{(9 m^2+35m+32)\;m\;m!}{24(m+4)!}\;.\quad\label{b6b}
\eea
In section \ref{r521}, formula (\ref{r521e}), we show that $L_r$ in (\ref{b6b}) 
come as coefficients in the series expansion of the function involving 
logarithmic and hypergeometric functions.
\begin{remark}\label{rem2}Consider an integer $N$ in the range $p^{\bar{r}}\leq 
N< p^{\bar{r}+1}$ where $\bar{r}=\left\lfloor\log_p N\right\rfloor$ and write 
(\ref{b4}) as follows,
\bea
F\left\{f;N,p^m\right\}=\sum_{r=0}^{\bar{r}-1}L_r\left(f;p^m\right)F\{f;N_r,1\}+
L_{\bar{r}}\left(f;p^m\right)F\{f;N_{\bar{r}},1\}\;.\label{b6c}
\eea
where $1\leq N_{\bar{r}}=\left\lfloor N/p^{\bar{r}}\right\rfloor<p$ and $F[f;
N_{\bar{r}},1]=\sum_{k=1}^{N_{\bar{r}}}f(k)<\sum_{k=1}^pf(k)$ is a finite 
number. 
\end{remark}
We make use of representation (\ref{b6c}) in section \ref{r21} when studying 
the asymptotics of renormalization functions for universality classes with 
$b_1<0$ (Lemma \ref{lem2}) and $b_1<b_2$ (Lemma \ref{lem4}).
\section{Renormalization Function with Simple Scaling}\label{r2}
Define the renormalization function as a ratio of two summatory functions
\bea
R\left(f;N,p^m\right)=\frac{F\left\{f;N,p^m\right\}}{F\{f;N,1\}}\;,\quad 
R\left(f;N,1\right)=1\;.\label{b7}
\eea
Substituting (\ref{i4}) into (\ref{b7}) we get its asymptotic behavior
\bea
R\left(f;N,p^m\right)\stackrel{N\to\infty}{=}\frac{{\cal R}_1\left(f;N,p^m
\right)}{1+{\cal O}\left(G_2(N)\right)}+\frac{{\cal R}_2\left(f;N,p^m\right)}
{1+{\cal O}\left(G_2(N)\right)}\;,\nonumber
\eea
where ${\cal R}_1\left(f;N,p^m\right)$ and ${\cal R}_2\left(f;N,p^m\right)$ are 
defined according to (\ref{b4}) as follows
\bea
{\cal R}_1\left(f;N,p^m\right)&=&\sum_{r=0}^{\left\lfloor\log_p N\right\rfloor}
L_r\left(f;p^m\right)\frac{G_1(N_r)}{G_1(N)}\;,\label{b9}\\
{\cal R}_2\left(f;N,p^m\right)&=&\sum_{r=0}^{\left\lfloor\log_p N\right\rfloor}
L_r\left(f;p^m\right)\frac{G_1(N_r)}{G_1(N)}\;{\cal O}\left(G_2(N_r)\right)\;.
\quad\label{b9b}
\eea
If both numerical series ${\cal R}_1\left(f;N,p^m\right)$ and ${\cal R}_2\left(
f;N,p^m\right)$ converge when $N\to\infty$, then
\bea
\lim_{N\to\infty}R\left(f;N,p^m\right):=R_{\infty}(f;p^m)=\lim_{N\to\infty}
{\cal R}_1\left(f;N,p^m\right)+\lim_{N\to\infty}{\cal R}_2\left(f;N,p^m\right)
\;.\label{b10}
\eea
What can be said about convergence of $R\left(f;N,p^m\right)$ without knowing
exactly the multiplicative function $f(k)$ itself ? Formulas (\ref{b9}) and 
(\ref{b9b}) for ${\cal R}_1\left(f;N,p^m\right)$ and ${\cal R}_2\left(f;N,p^m
\right)$ indicate that a large portion of information is hidden in the 
asymptotics $G_1(N)$ and $G_2(N)$. 

Substitute $G_1(N)=N^{a_1}\left(\log N\right)^{b_1}$ into (\ref{b9}) and get
\bea
{\cal R}_1\left(f;N,p^m\right)=\sum_{r=0}^{\left\lfloor\log_p N\right\rfloor}
\frac{L_r\left(f;p^m\right)}{p^{a_1r}}\left(1-\frac{r}{\log_pN}\right)^{b_1}.  
\label{b10a}
\eea
Regarding ${\cal R}_2\left(f;N,p^m\right)$, which is responsible for 
contribution of the error term into $R\left(f;N,p^m\right)$, note that according
to the definition (\ref{i5}) we get ${\cal O}\left(G_2(N)\right)\!\leq{\cal C}\;
G_2(N)$. Applying this inequality to formula (\ref{b9b}),
\bea
\left|{\cal R}_2\left(f;N,p^m\right)\right|\leq{\cal C}\sum_{r=0}^{\left\lfloor
\log_p N\right\rfloor}\left|L_r\left(f;p^m\right)\right|\;\frac{G_1(N_r)}
{G_1(N)}\;G_2(N_r)\;,\nonumber
\eea
and substituting $G_2(N)=N^{-a_2}\left(\log N\right)^{-b_2}$ into the last 
expression we get an estimate,
\bea
\left|{\cal R}_2\left(f;N,p^m\right)\right|\leq\frac{{\cal C}}{N^{a_2}\;\left(
\log N\right)^{b_2}}\;{\cal R}_3\left(f;N,p^m\right)\;,\quad\mbox{where}
\label{b10w}
\eea
\bea
{\cal R}_3\left(f;N,p^m\right)=\sum_{r=0}^{\left\lfloor\log_pN\right\rfloor}\!
\frac{\left|L_r\left(f;p^m\right)\right|}{p^{r(a_1-a_2)}}\left(1-\frac{r}{\log_p
N}\right)^{b_1-b_2}.\label{b10x}
\eea
In sections \ref{r21} and \ref{r22} we give a detailed analysis of convergence 
of ${\cal R}_1\left(f;N,p^m\right)$ and ${\cal R}_2\left(f;N,p^m\right)$ for 
multiplicative functions $f(k)$ of several universality classes. Start with a 
specific class of $f(k)$ and, assuming only $L_r\left(f;p^m\right)\geq 0$, prove
that the convergence of ${\cal R}_1\left(f;N,p^m\right)$ implies the convergence
of ${\cal R}_2\left(f;N,p^m\right)$ to zero.
\begin{proposition}\label{pro1}Let $f(k)\in{\mathbb B}\left\{N^{a_1}\left(\log 
N\right)^{b_1};\:N^{-a_1}\left(\log N\right)^{-b_2}\right\}$ be given such that 
$L_r\left(f;p^m\right)\geq 0$ and let ${\cal R}_1\left(f;N,p^m\right)$ be 
convergent. Then ${\cal R}_2\left(f;N,p^m\right)\stackrel{N\to\infty}
{\longrightarrow}0$, in each of the cases,
\bea
&&1)\quad a_1=a_2=0,\;b_1\geq b_2>0,\quad 2)\quad 0<a_2<a_1,\;0\leq b_2\leq b_1
\;\;\mbox{or}\;\;0=a_2<a_1,\;0<b_2\leq b_1,\nonumber\\
&&3)\quad 0<a_2<a_1,\;b_1<b_2,\;b_1<0\quad\mbox{or}\quad 0=a_2<a_1,\;b_1<0<b_2
\;.\nonumber
\eea
\end{proposition}
{\sf Proof} $\;\;\;$Keeping in mind $L_r\left(f;p^m\right)\geq 0$ and comparing 
(\ref{b10a}) and (\ref{b10x}) we conclude that if ${\cal R}_1\left(f;N,p^m
\right)$ is convergent when $a_1=0$ and $b_1\geq 0$ then ${\cal R}_3\left(f;N,
p^m\right)$ is also convergent when $a_1=a_2=0$ and $b_1\geq b_2$. Substituting
this into (\ref{b10w}) we get
\bea
\left|{\cal R}_2\left(f;N,p^m\right)\right|\leq {\cal C}\;\left(\log N\right)^{
-b_2}{\cal R}_3\left(f;N,p^m\right)\;,\label{b10w1}
\eea
that proves Proposition if $b_2>0$. Indeed, if $b_2$ is positive and ${\cal R}_
3\left(f;N,p^m\right)$ is convergent, then the right hand side (r.h.s.) in 
(\ref{b10w1}) is convergent to zero and so does ${\cal R}_2\left(f;N,p^m\right)$.
Applying similar arguments in the two others cases (2) and (3) we prove 
Proposition completely.$\;\;\Box$

Table 1 shows when the 1st item in Proposition \ref{pro1} can be applied: this
is the inverse Dedekind function: $a_1=a_2=0$, $b_1=b_2=1$ and by (\ref{b6a}) 
$L_r\left(1/\psi;p^m\right)\geq 0$. But neither the Euler totient function nor 
its inverse can be studied by Proposition \ref{pro1} which is quite weak 
statement.

The convergence of ${\cal R}_1\left(f;N,p^m\right)$ implies a zero limit of 
${\cal R}_2\left(f;N,p^m\right)$, $N\to\infty$, in much wider range of varying 
degrees $a_1,b_1$ and $a_2,b_2$. Indeed, to provide the convergence of ${\cal R}
_2\left(f;N,p^m\right)$ to zero there is no need to require the convergence of 
${\cal R}_3\left(f;N,p^m\right)$ in (\ref{b10w}) but rather to allow a growth of
${\cal R}_3\left(f;N,p^m\right)$ with a rate less than $N^{a_2}\left(\log N
\right)^{b_2}$. However, this would require more assumptions about $L_r\left(f;
p^m\right)$. In the next sections \ref{r21} and \ref{r22} we study the 
convergence problem in more details and prove the main result of this section in
Theorem \ref{thm1}.
\subsection{Convergence of ${\cal R}_1\left(f;N,p^m\right)$}\label{r21}
In this section we study the convergence of ${\cal R}_1\left(f;N,p^m\right)$ 
when $N\to\infty$. The conditions imposed on $L_r\left(f;p^m\right)$ provide 
convergence of ${\cal R}_1\left(f;N,p^m\right)$. Throughout this and the next 
sections we repeatedly make use of the squeeze ({\sf SQ}) theorem \cite{spi80},
which is also known as the pinching or sandwich theorems.
\begin{lemma}\label{lem1}Let a function $f(k)\in{\mathbb B}\left\{N^{a_1}\left(
\log N\right)^{b_1};\;G_2(N)\right\}$, $a_1\geq 0$, $b_1\geq 0$, be given and
let there exist two numbers ${\cal K}>0$ and $\gamma<a_1$ and an integer $r_*
\geq 0$ such that $|L_r\left(f;p^m\right)|\leq {\cal K}p^{\gamma r}$ for all 
$r\geq r_*$. Then 
\bea
\lim_{N\to\infty}{\cal R}_1\left(f;N,p^m\right)=\sum_{r=0}^{\infty}\frac{L_r
\left(f;p^m\right)}{p^{a_1r}}\;.\label{b11}
\eea
\end{lemma}
{\sf Proof} $\;\;\;$First, consider ${\cal R}_1\left(f;N,p^m\right)$ given in 
(\ref{b10a}) when $a_1\geq 0$, $b_1\in{\mathbb Z}_+\cup\{0\}$, i.e., $b_1$ is
a nonnegative integer. After binomial expansion in (\ref{b10a}) we get,
\bea
{\cal R}_1\left(f;N,p^m\right)=\sum_{r=0}^{\left\lfloor\log_p N\right\rfloor}  
\frac{L_r\left(f;p^m\right)}{p^{a_1r}}+\sum_{k=1}^{b_1}\left(\frac{-1}{\log_pN} 
\right)^k{b_1\choose k}\sum_{r=0}^{\left\lfloor\log_p N\right\rfloor}\frac{
L_r\left(f;p^m\right)}{p^{a_1r}}\;r^k\;.\label{b10b}
\eea
Focus on the inner sum
\bea
{\cal R}_4\left(f;N,p^m,k\right)=\frac1{\left(\log_pN\right)^k}\sum_{r=0}^{\left
\lfloor\log_p N\right\rfloor}\frac{L_r\left(f;p^m\right)}{p^{a_1r}}\;r^k,\quad 
1\leq k\leq b_1\;,\label{b10c}
\eea
and prove that the sum in (\ref{b10c}) is convergent absolutely. To find an 
estimate for ${\cal R}_4\left(f;N,p^m,k\right)$ we have to consider the last 
sum at interval $\left(r_*,\left\lfloor\log_p N\right\rfloor\right)$ where an 
inequality $|L_r\left(f;p^m\right)|\leq {\cal K}p^{\gamma r}$ holds. However, 
because of the prefactor $\left(\log_pN\right)^{-k}$, $k\geq 1$, a summation at 
interval $(0,r_*-1)$ does not contribute to the limit when $N\to\infty$ and does
not change the convergence of the entire sum (\ref{b10c}). Then
\bea
|{\cal R}_4\left(f;N,p^m,k\right)|\leq\frac1{\left(\log_pN\right)^k}\sum_{r=0}^{
\left\lfloor\log_pN\right\rfloor}\frac{|L_r\left(f;p^m\right)|}{p^{a_1r}}\;r^k
\leq\frac{{\cal K}}{\left(\log_pN\right)^k}\!\!\sum_{r=0}^{\left\lfloor\log_pN
\right\rfloor}\!\!\frac{r^k}{p^{\epsilon r}}\;,\label{b11v}
\eea
where $\epsilon=a_1-\gamma>0$. Denote $M=\left\lfloor\log_p N\right\rfloor$ and 
consider the sum in the r.h.s. of (\ref{b11v}),
\bea
T(p,k,\epsilon,M)=\sum_{r=0}^M\frac{r^k}{p^{\epsilon r}}={\rm Li}_{-k}\left(
p^{-\epsilon}\right)-p^{-\epsilon(M+1)}\Phi\left(p^{-\epsilon}\!,-k,M+1\right)
\;,\label{b11c1}
\eea
where ${\rm Li}_s(z)$ and $\Phi(z,s,a)$ are the polylogarithm function and the
Hurwitz-Lerch zeta function \cite{em81}, respectively,
\bea
{\rm Li}_s(z)=\sum_{k=1}^{\infty}\frac{z^k}{k^s}\;,\hspace{1cm}
\Phi(z,s,a)=\sum_{k=0}^{\infty}\frac{z^k}{(a+k)^s}\;.\label{b11a1}
\eea
Keeping in mind the asymptotics of $\Phi(z,s,a)$, $a\to\infty$, for fixed $s$
and $z$ (Thm.1, \cite{fl04}), $\Phi(z,s,a)\simeq a^{-s}/(1-z)$, and combining it
with (\ref{b10c}) and (\ref{b11c1}), we get
\bea
|{\cal R}_4\left(f;N,p^m,k\right)|\stackrel{N\to\infty}{\leq}{\cal K}\left[
\frac{{\rm Li}_{-k}\left(p^{-\epsilon}\right)}{\left(\log_pN\right)^k}-\frac{N^{
-\epsilon}}{p^{\epsilon}-1}\right],\quad 1\leq k\leq b_1\;.\label{b11b1}
\eea
Since the polylogarithm ${\rm Li}_{-k}\left(z\right)$, $0\leq k<\infty$, is a
bounded rational function in $z$ when $0\leq z<1$, then for $\epsilon>0$ the
function ${\rm Li}_{-k}\left(p^{-\epsilon}\right)$ is also bounded. Thus, ${\cal
R}_4\left(f;N,p^m,k\right)$ is convergent to zero and by (\ref{b10b}) and
(\ref{b10c}) the limit (\ref{b11}) holds.  

Extend this result on the entire set of the nonnegative real numbers $b_1$. In 
order to do this, we make use of the {\sf SQ} theorem and start with trivial 
inequalities,
\bea
\left(1-\frac{r}{\log_pN}\right)^{\left\lfloor b_1+1\right\rfloor}\leq\left(1-
\frac{r}{\log_pN}\right)^{b_1}\leq\left(1-\frac{r}{\log_pN}\right)^{\left\lfloor
b_1\right\rfloor},\quad 0\leq r\leq\left\lfloor\log_p N\right\rfloor\;,  
\label{a218}
\eea
which due to (\ref{b10a}) implies the following relations,
\bea
&&{\cal J}_1\left(f;p^m,\log_pN,\left\lfloor b_1+1\right\rfloor\right)\leq 
{\cal  J}_1\left(f;p^m,\log_pN,b_1\right)\leq {\cal J}_1\left(f;p^m,\log_pN,
\left\lfloor b_1\right\rfloor\right)\;,\quad\label{a219a}\\
&&\mbox{where}\hspace{1cm}{\cal J}_1\left(f;p^m,M,b\right)=\sum_{r=0}^{\left
\lfloor M\right\rfloor}|L_r\left(f;p^m\right)|\;p^{-a_1r}\left(1-\frac{r}{M}
\right)^b\;,\quad b\geq 0\;.\nonumber
\eea
By proof on convergence of ${\cal R}_1\left(f;N,p^m\right)$ with nonnegative 
integer degrees $b_1$ and by (\ref{a219a}) and by the {\sf SQ} theorem it 
follows the convergence of ${\cal R}_1\left(f;N,p^m\right)$ with real 
$b_1\geq 0$.$\;\;\;\;\;\;\Box$

Lemma \ref{lem1} can be applied to the totient functions $\phi(k)$, $\psi(k)$ 
and their inverse $1/\phi(k)$, $1/\psi(k)$ with the functions $L_r\left(f;p^m
\right)$ calculated in (\ref{b6a}).
\begin{lemma}\label{lem2}Let a function $f(k)\in{\mathbb B}\left\{N^{a_1}\left(
\log N\right)^{b_1};\;G_2(N)\right\}$, $a_1>0$, $b_1<0$, be given and let there 
exist two numbers ${\cal K}>0$ and $\gamma<a_1$ and an integer $r_*\geq 0$ such 
that $|L_r\left(f;p^m\right)|\leq {\cal K}p^{\gamma r}$ for all $r\geq r_*$. 
Then (\ref{b11}) holds.
\end{lemma}
{\sf Proof} $\;\;\;$Consider ${\cal R}_1\left(f;N,p^m\right)$ given in  
(\ref{b10a}) when $a_1>0$, $b_1\in{\mathbb Z}_-$, i.e., $b_1$ is a negative 
integer. In order to avoid its divergence at $r=\log_p N$ we use the 
representation (\ref{b6c}) in Remark \ref{rem2},
\bea
{\cal R}_1\left(f;N,p^m\right)=\!\!\!\sum_{r=0}^{\left\lfloor\log_p N\right
\rfloor-1}\!\!\frac{L_r\left(f;p^m\right)}{p^{a_1r}}\left(1-\frac{r}{\log_pN}
\right)^{-|b_1|}\!\!+\frac{A}{{\cal F}}\frac{\left(\log_p N\right)^{|b_1|}}
{N^{a_1}}\;,\label{b20}
\eea
where $A=L_{\bar{r}}\left(f;p^m\right)F[f;N_{\bar{r}},1]<\infty$ and $\bar{r}=
\left\lfloor\log_p N\right\rfloor$ were defined in (\ref{b6c}).
 
Estimate ${\cal R}_1\left(f;N,p^m\right)$ when $|L_r\left(f;p^m\right)|<{\cal K}
p^{\gamma r}$, ${\cal K}>0$ and $\gamma<a_1$. Note that the last term in 
(\ref{b20}) does not contribute to the asymptotics of ${\cal R}_1\left(f;N,p^m
\right)$ when $N\to\infty$ and therefore it can be skipped hereafter. We use an 
identity 
\bea
(1-x)^{-b}=1+x\sum_{k=1}^b(1-x)^{-k}\;,\quad b\in{\mathbb Z}_+\;,\label{b218a}
\eea
and represent (\ref{b20}) as follows,
\bea
{\cal R}_1\left(f;N,p^m\right)\simeq\sum_{r=0}^{\left\lfloor\log_p N\right
\rfloor-1}\frac{L_r\left(f;p^m\right)}{p^{a_1r}}+\sum_{k=1}^{|b_1|}{\cal R}_5
\left(f;N,p^m,k\right),\hspace{1cm}\mbox{where}\label{b20a}
\eea
\bea
{\cal R}_5\left(f;N,p^m,k\right)=\frac1{\log_p N}\sum_{r=0}^{\left\lfloor\log_p
N\right\rfloor-1}\frac{L_r\left(f;p^m\right)}{p^{a_1r}}\;r\left(1-\frac{r}
{\log_pN}\right)^{-k}.\label{b20b}
\eea
Note that the following inequality holds,
\bea
\left(1-\frac{r}{\log_pN}\right)^{-1}\leq 1+r\,\quad\mbox{if}\quad0\leq r<
\left\lfloor\log_p N\right\rfloor\;.\label{b20i}
\eea
Substitute (\ref{b20i}) into (\ref{b20b}) and keep in mind that due to the
prefactor $\left(\log_pN\right)^{-1}$ in the r.h.s. of Eq. (\ref{b20b}) the same
convergence of ${\cal R}_5\left(f;N,p^m,k\right)$ holds at intervals $\left(r_*,
\left\lfloor\log_p N\right\rfloor-1\right)$ and $\left(0,\left\lfloor\log_p N
\right\rfloor-1\right)$ (see discussion in proof of Lemma \ref{lem1}). Then we 
arrive at estimate,
\bea
|{\cal R}_5\left(f;N,p^m,k\right)|\leq\frac{{\cal K}}{\log_p N}\!\!\sum_{r=0}^{
\left\lfloor\log_p N\right\rfloor-1}\!\!\frac{r(1+r)^k}{p^{\epsilon r}}=
\frac{{\cal K}}{\log_p N}\sum_{j=0}^k\;{k\choose j}\!\sum_{r=0}^{\left\lfloor
\log_p N\right\rfloor-1}\!\!\frac{r^{j+1}}{p^{\epsilon r}}\;.\nonumber
\eea
The rest of the proof follows by applying the same arguments of asymptotics of 
the Hurwitz-Lerch zeta function, as it was done in Lemma \ref{lem1}, 
\bea
|{\cal R}_5\left(f;N,p^m,k\right)|\leq {\cal K}\sum_{j=0}^k{k\choose j}\left[
\frac{{\rm Li}_{-(j+1)}\left(p^{-\epsilon}\right)}{\log_pN}-\frac1{p^{\epsilon}
-1}\frac{\left(\log_pN\right)^{j}}{N^{\epsilon}}\right]\;.\label{b20g}
\eea
By comparison the r.h.s. in (\ref{b20g}) and (\ref{b11b1}) we conclude that 
${\cal R}_5\left(f;N,p^m,k\right)$ is convergent to zero when $N\to\infty$.
Thus, by (\ref{b20a}) the limit (\ref{b11}) holds for $b_1\in{\mathbb Z}_-$.

We extend this result by the {\sf SQ} theorem on all negative real $b_1$. This
can be done by inequality (\ref{a219a}) for another function ${\cal J}_2\left(
f;p^m,M,b\right)$,
\bea
&&{\cal J}_2\left(f;p^m,\log_pN,\left\lfloor b_1+1\right\rfloor\right)\leq
{\cal J}_2\left(f;p^m,\log_pN,b_1\right)\leq {\cal J}_2\left(f;p^m,\log_pN,
\left\lfloor b_1\right\rfloor\right)\;,\label{b20j}\\
&&\mbox{where}\hspace{1cm}
{\cal J}_2\left(f;p^m,M,b\right)=\sum_{r=0}^{\left\lfloor M\right\rfloor-1}
\frac{|L_r\left(f;p^m\right)|}{p^{a_1r}}\left(1-\frac{r}{M}\right)^b\;,\quad 
b<0\;.\nonumber
\eea
By proof on convergence of ${\cal R}_1\left(f;N,p^m\right)$ with negative 
integer degrees $b_1$ and by (\ref{b20j}) and by the {\sf SQ} theorem it 
follows the convergence of ${\cal R}_1\left(f;N,p^m\right)$ with real $b_1< 0$.
$\;\;\;\;\;\;\Box$
\subsection{Convergence of ${\cal R}_2\left(f;N,p^m\right)$}\label{r22}
In this section we consider the convergence of ${\cal R}_2\left(f;N,p^m\right)$,
$N\to\infty$, defined in (\ref{b9}) and responsible for contribution of the 
error term to $R_{\infty}\left(f;p^m\right)$. 
\begin{lemma}\label{lem3}Let a function $f(k)\in{\mathbb B}\left\{N^{a_1}\left(
\log N\right)^{b_1};N^{-a_2}\left(\log N\right)^{-b_2}\right\}$, $b_1\geq b_2$,
be given and let there exist two numbers ${\cal K}>0$ and $\gamma<a_1$ and an 
integer $r_*\geq 0$ such that $|L_r\left(f;p^m\right)|\leq {\cal K}p^{\gamma r}$
for all $r\geq r_*$, then
\bea
\lim_{N\to\infty}{\cal R}_2\left(f;N,p^m\right)=0\;.\label{q0}
\eea
\end{lemma}
{\sf Proof} $\;\;\;$Denote $b_1-b_2=e$ and make use of a simple inequality, 
$(1-x)^e\leq 1$ when $0\leq x\leq 1$, $e\geq 0$. Formula (\ref{b10w}) can be
rewritten as follows,
\bea
\left|{\cal R}_2\left(f;N,p^m\right)\right|\leq\frac{{\cal C}}{N^{a_2}\left(
\log N\right)^{b_2}}\sum_{r=0}^{\left\lfloor\log_pN\right\rfloor}\frac{\left|
L_r\left(f;p^m\right)\right|}{p^{(a_1-a_2)r}}\leq\frac{{\cal C}{\cal K}}{N^{a_2}
\left(\log N\right)^{b_2}}\sum_{r=0}^{\left\lfloor\log_pN\right\rfloor}p^{\nu
r}\;,\label{q7}
\eea
where $\nu=a_2-a_1+\gamma$. Keep in mind that due to the prefactor $N^{-a_2}
\left(\log N\right)^{-b_2}$ in (\ref{q7}), the same convergence of the r.h.s. 
in (\ref{q7}) holds at intervals $\left(r_*,\left\lfloor\log_p N\right\rfloor
\right)$ and $\left(0,\left\lfloor\log_p N\right\rfloor\right)$ (see discussion 
in proof of Lemma \ref{lem1}). 

The further calculations are dependent on the sign of $\nu$. If $\nu<0$ then
\bea
\gamma<a_1-a_2\;,\quad\left|{\cal R}_2\left(f;N,p^m\right)\right|\leq\frac{{\cal
C}{\cal K}}{1-p^{\nu}}\frac{N^{-a_2}}{\left(\log N\right)^{b_2}}\;.\label{q8}
\eea
If $\nu=0$ then
\bea
\gamma=a_1-a_2\;,\quad\left|{\cal R}_2\left(f;N,p^m\right)\right|\leq {\cal C}
{\cal K}\frac{N^{\gamma-a_1}}{\left(\log N\right)^{b_2-1}}\;.\label{q9}
\eea
Finally, if $\nu>0$ then
\bea
\gamma>a_1-a_2\;,\quad\left|{\cal R}_2\left(f;N,p^m\right)\right|\leq\frac{{\cal
C}{\cal K}}{N^{a_2}\left(\log N\right)^{b_2}}\frac{p^{\nu\left(\left\lfloor
\log_pN\right\rfloor+1\right)}}{p^{\nu}-1}\simeq\frac{{\cal C}{\cal K}}{1-p^{-
\nu}}\frac{N^{\gamma-a_1}}{\left(\log N\right)^{b_2}}\;.\quad\label{q10}
\eea

Require now that all r.h.s. in (\ref{q8}), (\ref{q9}) and (\ref{q10}) converge 
to zero when $N\to\infty$. Regarding the 1st case (\ref{q8}) this always holds 
because by (\ref{i7}) if $a_2=0$ then $b_2>0$, and if $a_2>0$ then $b_2\geq 0$, 
so the r.h.s. in (\ref{q8}) is decreasing function. So, it results in 
requirement, $\gamma<a_1-a_2$. In two other cases we have necessary conditions,
\bea
&&a_1-a_2=\gamma\leq a_1\quad\mbox{if}\quad b_2>1\;,\quad\mbox{and}\quad
a_1-a_2=\gamma<a_1\quad\mbox{if}\quad-\infty<b_2<\infty\;,\quad\label{q11}\\
&&a_1-a_2<\gamma\leq a_1\quad\mbox{if}\quad b_2>0\;,\quad\mbox{and}\quad
a_1-a_2<\gamma<a_1\quad\mbox{if}\quad-\infty<b_2<\infty\;.\quad\label{q12}
\eea
Summarizing the necessary conditions (\ref{q11}), (\ref{q12}) and $\gamma<a_1-
a_2$, we conclude that the numerical series ${\cal R}_2\left(f;N,p^m\right)$ is 
convergent to zero absolutely when $N\to\infty$ and irrespectively to the sign 
of $b_2$ if $\gamma<a_1$. This proves formula (\ref{q0}).$\;\;\;\;\;\;\Box$
\begin{lemma}\label{lem4}Let a function $f(k)\in{\mathbb B}\left\{N^{a_1}\left(
\log N\right)^{b_1};N^{-a_2}\left(\log N\right)^{-b_2}\right\}$, $b_1<b_2$, be
given and let there exist two numbers ${\cal K}>0$ and $\gamma<a_1$ and an 
integer $r_*\geq 0$ such that $|L_r\left(f;p^m\right)|\leq {\cal K}p^{\gamma r}$
for all $r\geq r_*$, then (\ref{q0}) holds.
\end{lemma}
{\sf Proof} $\;\;\;$Consider ${\cal R}_3\left(f;N,p^m\right)$ given in 
(\ref{b10x}) and, according to Remark \ref{rem2}, rewrite it as follows,
\bea
{\cal R}_3\left(f;N,p^m\right)=\sum_{r=0}^{\left\lfloor\log_p N\right\rfloor-1}
\frac{\left|L_r\left(f;p^m\right)\right|}{p^{(a_1-a_2)r}}\left(1-\frac{r}{\log_p
N}\right)^{-|b_2-b_1|}+\frac{A}{{\cal F}}\frac{N^{-a_1}}{\left(\log_p N\right)^{
b_1}}\;,\label{q13}
\eea
where the upper bound in the sum is taken in order to avoid its divergence and 
$A$ is defined in (\ref{b20}), more details see in proof of Lemma \ref{lem2}, 
formula (\ref{b20}). The last term in (\ref{q13}) does not contribute to 
asymptotics of ${\cal R}_2\left(f;N,p^m\right)$ when $N\to\infty$ and can be 
skipped. Indeed, if $a_1\geq 0$, $b_1>0$ it converges to zero when $N\to\infty$;
in the case $a_1=b_1=0$ acording to (\ref{i7}) we have $a_2>0$ or $a_2=0$, $b_2>
0$ that again makes it irrelevant due to prefactor $N^{-a_2}\left(\log_pN\right)
^{-b_2}$ in formula (\ref{b10w}) for ${\cal R}_2\left(f;N,p^m\right)$. Apply an 
inequality (\ref{b20i}) in the range $0\leq r\leq\left\lfloor\log_p N\right
\rfloor-1$,
\bea
\left(1-\frac{r}{\log_p N}\right)^{-(b_2-b_1)}\leq\left(\log_p N\right)^{b_2-
b_1}\;,\nonumber
\eea
and substitute it into (\ref{q13})
\bea
{\cal R}_3\left(f;N,p^m\right)\leq\left(\log_p N\right)^{b_2-b_1}\sum_{r=0}^{
\left\lfloor\log_p N\right\rfloor-1}\frac{\left|L_r\left(f;p^m\right)\right|}{
p^{r(a_1-a_2)}}\;.\label{q14} 
\eea
If $a_2>0$ we apply to (\ref{q14}) the constraints on $L_{r}\left(f;p^{m}
\right)$ and substitute the result into (\ref{q69}),
\bea
{\cal R}_2\left(f;N,p^m\right)\leq\frac{{\cal C}{\cal K}}{N^{a_2}\;\left(\log_p
N\right)^{b_1}}\sum_{r=0}^{\left\lfloor\log_p N\right\rfloor-1}p^{\nu r}\;,\quad
\nu=a_2-a_1+\gamma\;.\label{q15}
\eea
By comparison (\ref{q15}) with (\ref{q7}) from Lemma \ref{lem3} we obtain
according to (\ref{q8}), (\ref{q9}) and (\ref{q10})
\bea
\left|{\cal R}_2\left(f;N,p^m\right)\right|&\leq&\frac{{\cal C}{\cal K}}{p^{\nu}
-1}\frac{N^{\gamma-a_1}}{\left(\log_c N\right)^{b_1}}\;,\quad \nu>0\;,
\label{q15a}\\
\left|{\cal R}_2\left(f;N,p^m\right)\right|&\leq&{\cal C}{\cal K}\;\frac{N^{-
a_2}}{\left(\log_c N\right)^{b_1-1}}\;,\quad \nu=0\;.\label{q15b}\\
\left|{\cal R}_2\left(f;N,p^m\right)\right|&\leq&\frac{{\cal C}{\cal K}}{1-p^{
\nu}}\;\frac{N^{-a_2}}{\left(\log_c N\right)^{b_1}}\;,\quad \nu<0\;.\label{q15c}
\eea
Thus, by (\ref{q15a}), (\ref{q15b}), (\ref{q15c}) and inequalities $\gamma<a_1$,
$a_2>0$ a series ${\cal R}_2\left(f;N,p^m\right)$ is convergent to zero when 
$N\to\infty$ irrespectively to the value of $b_1$.

Consider another case $a_2=0$ and compare formulas (\ref{b20}) and (\ref{q13}) 
for $b_1<0$ and $b_1<b_2$, respectively. A difference in degrees, $b_1-b_2$ and 
$b_1$, does not break the main result of Lemma \ref{lem2}: only the 1st leading 
term in (\ref{b20a}) is survived when $N\to\infty$. When we apply it to 
(\ref{q13}) and make use of constraint on $L_{r}\left(f;p^m\right)$ we get,
\bea
{\cal R}_3\left(f;N,p^m\right)\stackrel{N\to\infty}{\simeq}\sum_{r=0}^{\left
\lfloor\log_p N\right\rfloor-1}\frac{L_r\left(f;p^m\right)}{p^{a_1r}}\leq{\cal K
}\sum_{r=0}^{\infty}p^{-\epsilon r}=\frac{{\cal K}}{1-p^{-\epsilon}}\;,\nonumber
\eea
where $\epsilon=a_1-\gamma>0$. Substituting the last estimate into (\ref{q69}) 
we obtain,
\bea
\left|{\cal R}_2\left(f;N,p^m\right)\right|\leq\frac{{\cal C}{\cal K}}{1-p^{-
\epsilon}}\cdot\left(\log_p N\right)^{-b_2}\;.\label{q16}
\eea
Recall that by (\ref{i7}) the degrees of the error term satisfy: if $a_2=0$ than
$b_2>0$. Thus, by (\ref{q16}) the series ${\cal R}_2\left(f;N,p^m\right)$ is 
convergent to zero when $N\to\infty$ that proves Lemma.$\;\;\;\;\;\;\Box$

Summarize the results of four Lemmas \ref{lem1}, \ref{lem2}, \ref{lem3} and 
\ref{lem4} on convergence of the numerical series ${\cal R}_1\left(f;N,p^m
\right)$ and ${\cal R}_2\left(f;N,p^m\right)$.
\begin{theorem}\label{thm1}Let a function $f(k)\in{\mathbb B}\left\{G_1(N);G_2
(N)\right\}$ be given with $a_i$ and $b_i$ satisfying (\ref{i7}) and let there 
exist two numbers ${\cal K}>0$ and $\gamma<a_1$ and an integer $r_*\geq 0$ such 
that $|L_r\left(f;p^m\right)|\leq {\cal K}p^{\gamma r}$ for all $r\geq r_*$. 
Then
\bea
R_{\infty}\left(f,p^m\right)=\sum_{r=0}^{\infty}\frac{L_r\left(f;p^m\right)}
{p^{a_1r}}\;.\label{q21}
\eea
\end{theorem}
{\sf Proof} $\;\;\;$According to Lemmas \ref{lem1} and \ref{lem2} if there exist
two numbers ${\cal K}>0$ and $\gamma<a_1$ and an integer $r_*\geq 0$ such that 
$|L_r\left(f;p^m\right)|\leq{\cal K}p^{\gamma r}$ for all $r\geq r_*$ then 
${\cal R}_1\left(f;N,p^m\right)$ is convergent to $\sum_{r=0}^{\infty}L_r\left(
f;p^m\right)p^{-a_1r}$ in the whole range (\ref{i7}) of varying parameters 
$a_1,b_1$ and $a_2,b_2$. On the other hand, according to Lemmas \ref{lem3} and 
\ref{lem4} by the same sufficient conditions the numerical series ${\cal R}_2
\left(f;N,p^m\right)$ is convergent to zero in the same range (\ref{i7}) of 
varying parameters. Then, in accordance with (\ref{b10}) we arrive at 
(\ref{q21}).$\;\;\;\;\;\;\Box$

There are a few questions on convergence of ${\cal R}_1\left(f;N,p^m\right)$ and
${\cal R}_2\left(f;N,p^m\right)$ which have been left open beyond the scope of 
Theorem \ref{thm1}. First, this is a problem of necessary convergence conditions
which need further discussion. Another question arises in view of convergence 
to zero of ${\cal R}_2\left(f;N,p^m\right)$: it can happen that for some $f(k)$ 
both renormalization functions ${\cal R}_1\left(f;N,p^m\right)$ and ${\cal R}_2
\left(f;N,p^m\right)$ are convergent to the non zero values while $L_r\left(f;
p^m\right)$ is satisfying less strong conditions than those given in Theorem 
\ref{thm1}. Thus, the following question has been left still open,
\begin{question}\label{que1}
Does the error term contribute to the renormalization function $R_{\infty}
\left(f;p^m\right)$ and what multiplicative arithmetic functions $f(k)$ can 
provide an affirmative answer ?
\end{question}
\subsection{Rational Representation of Renormalization Function}\label{r23}
In sections \ref{r21} and \ref{r22} we have found the requirements which suffice
to make ${\cal R}_1\left(f;N,p^m\right)$ convergent and ${\cal R}_2\left(f;N,p^m
\right)$ vanishing. These conditions are presented through the characteristic 
functions $L_r\left(f;p^m\right)$ given recursively in (\ref{b5}). Their 
straightforward formulas (\ref{b6}) look cumbersome and lead in particular cases
to rather complicate expressions, e.g., (\ref{b6b}). This is why in this section
we give another representation for $R_{\infty}\left(f;p^m\right)$ avoiding the  
use of $L_r\left(f;p^m\right)$.

Substituting (\ref{b5}) into (\ref{q21}) we get an infinite series for 
$R_{\infty}\left(f;p^m\right)\equiv R_{\infty}$,
\bea
R_{\infty}=f\left(p^m\right)+\frac{f\left(p^{m+1}\right)-
L_0\left(f;p^m\right)f(p)}{p^{a_1}}+\frac{f\left(p^{m+2}\right)-L_1\left(f;p^m
\right)f(p)-L_0\left(f;p^m\right)f(p^2)}{p^{2a_1}}+\nonumber\\
\frac{f\left(p^{m+3}\right)-L_2\left(f;p^m\right)f(p)-L_1\left(f;p^m\right)
f(p^2)-L_0\left(f;p^m\right)f(p^3)}{p^{2a_1}}+\ldots\;.\nonumber
\eea
Recasting the terms in the last expression we obtain
\bea
R_{\infty}\left(f;p^m\right)=\sum_{r=0}^{\infty}\frac{f\left(p^{m+r}\right)}
{p^{a_1r}}-\frac{f(p)}{p^{a_1}}\sum_{r=0}^{\infty}\frac{L_r\left(f;p^m\right)}
{p^{a_1r}}-\frac{f(p^2)}{p^{2a_1}}\sum_{r=0}^{\infty}\frac{L_r\left(f;p^m
\right)}{p^{a_1r}}-\ldots\;.\label{q22}
\eea
Introduce two numerical series ${\bf U}\left(f;p^m\right)$ and ${\bf V}\left(
f;p\right)$
\bea
{\bf U}\left(f;p^m\right)=\sum_{r=0}^{\infty}\frac{f\left(p^{r+m}\right)}{p^{
a_1r}}\;,\hspace{1cm}{\bf V}\left(f;p\right)=\sum_{r=1}^{\infty}\frac{f\left(
p^{r}\right)}{p^{a_1r}}\;,\label{q23}
\eea
and assume that they are convergent. Then, by comparison of (\ref{q22}) and 
(\ref{q21}) we get
\bea
R_{\infty}\left(f;p^m\right)=\frac{{\bf U}\left(f;p^m\right)}{1+{\bf V}\left(f;
p\right)}\;.\label{q24}
\eea
if the denominator in (\ref{q24}) does not vanish. For short it can be written 
as follows, ${\bf V}\left(f;p\right)+1=\sum_{r=0}^{\infty}f\left(p^{r}\right)
p^{-a_1r}$. However, we prefer to stay with (\ref{q24}), otherwise one can make
an error in calculations, e.g., $\phi\left(p^r\right)=p^{r-1}(p-1)$, $r\geq 
1$, and $\phi\left(p^0\right)=1$, but $\phi\left(p^0\right)\neq p^{-1}(p-1)$.

Formula (\ref{q24}) gives a rational representation of the renormalization 
function $R_{\infty}\left(f;p^m\right)$ which is free of intermediate 
calculations of $L_r\left(f;p^m\right)$. What can be said about convergence of 
$R_{\infty}\left(f;p^m\right)$ in terms of $f\left(p^r\right)$? 
\begin{theorem}\label{thm2}Let a function $f(k)\in{\mathbb B}\left\{G_1(N);G_2
(N)\right\}$ be given and let there exist two numbers ${\cal K}_1>0$ and
$\gamma_1<a_1$ and an integer $r_*\geq 0$ such that $|f\left(p^r\right)|\leq
{\cal K}_1p^{\gamma_1 r}$ for all $r\geq r_*$. Then $R_{\infty}\left(f;p^m
\right)$ is convergent in accordance with (\ref{q24}) if ${\bf V}\left(f;p
\right)+1\neq 0$.
\end{theorem}
{\sf Proof} $\;\;\;$The constraints $|f\left(p^r\right)|\leq{\cal K}_1p^{\gamma
_1r}$, ${\cal K}_1>0$ and $\gamma_1<a_1$, for all $r\geq r_*$ implies the 
convergence of ${\bf U}\left(f;p^m\right)$ and ${\bf V}\left(f;p\right)$ that 
follows by their absolute convergence,
\bea
\left|{\bf U}\left(f;p^m\right)\right|\leq{\cal K}_1p^{\gamma_1m}\sum_{r=0}^{
\infty}p^{(\gamma_1-a_1)r}=\frac{{\cal K}_1p^{\gamma_1m}}{1-p^{\gamma_1-a_1}}\;,
\hspace{.6cm}
\left|{\bf V}\left(f;p\right)\right|\leq{\cal K}_1\sum_{r=0}^{\infty}p^{(
\gamma_1-a_1)r}=\frac{{\cal K}_1}{1-p^{\gamma_1-a_1}}\;.\nonumber
\eea
Thus, if the denominator in (\ref{q24}) does not vanish, ${\bf V}\left(f;p
\right)+1\neq 0$, then $R_{\infty}\left(f;p^m\right)$ is convergent in 
accordance with (\ref{q24}).$\;\;\;\;\;\;\Box$

In Table 3 we present different multiplicative functions and their corresponding
parameters $a_1$, $\gamma_1$, ${\cal K}_1$ and $r_*$. All functions satisfy the 
constraints of Theorem \ref{thm2} for $p\geq 2$.
\begin{center}
{\bf Table$\;$3}. Multiplicative functions and their corresponding parameters 
$a_1$, $\gamma_1$, ${\cal K}_1$ and $r_*$.\\
\vspace{.5cm}
\begin{tabular}{|c||c|c|c|c|c|c|c|c|} \hline
$f(k)$ & $\phi(k)$ & $1/\phi(k)$ & $\psi(k)$ & $1/\psi(k)$ & $J_n(k)$ & 
$\mu^2(k)$ & $2^{\omega(k)}$ & $3^{\omega(k)}$ \\\hline\hline
$a_1$ & 2 & 0 & 2 & 0 & $n+1$ & 1 & 1 & 1\\\hline
$\gamma_1$ & 1 & -1 & 1 & -1 & $n$ & 0 & 0 & 0\\\hline
${\cal K}_1$ & $(p-1)/p$ & $p/(p-1)$ & $(p+1)/p$ & $p/(p+1)$ & $(p^n-1)/p^n$ & 
1 & 2 & 3 \\\hline $r_*$ & 1 & 1 & 1 & 1 & 1 & 0 & 0 & 0 \\\hline
\end{tabular}
\label{ta3}
\end{center}
Regarding the convergence of ${\bf U}\left(f;p^m\right)$ and ${\bf V}\left(f;p
\right)$ defined in (\ref{q23}) we present an example which shows that the 
renormalization function $R_{\infty}\left(f;p^m\right)$ can exist even when both
${\bf U}\left(f;p^m\right)$ and ${\bf V}\left(f;p\right)$ are divergent. 
Consider $f(k)=2^{\Omega(k)}$, $2^{\Omega(p^r)}=2^r$ and get
\bea
\left\{\begin{array}{l}a_1=\gamma_1=1\\{\cal K}_1=1\end{array}\right.\!\!,\;
\left\{\begin{array}{l}L_0\left(2^{\Omega(k)};p^m\right)=2^m\\L_r\left(2^{
\Omega(k)};p^m\right)=0,\;r\geq 1\end{array}\right.\!\!,\;
\left\{\begin{array}{l}{\bf U}\left(f;2^m\right)\stackrel{N\to\infty}{\simeq}
\infty\\
{\bf V}\left(f;2\right)\stackrel{N\to\infty}{\simeq}\infty\end{array}\right.
\!\!,\;R_{\infty}\left(2^{\Omega(k)},2^m\right)=2^m\;.\nonumber
\eea
Here the conditions of Theorem \ref{thm1} are satisfied for all $r\geq 1$ and 
only the 1st term is left nonzero in series (\ref{q21}). Both ${\bf U}\left(f;
2^m\right)$ and ${\bf V}\left(f;2\right)$ are divergent, e.g., ${\bf V}\left(f;2
\right)=\sum_{r=1}^{\infty}1$, and therefore Theorem \ref{thm2} cannot be 
applied. In other words, Theorem \ref{thm1} has much wider area of application 
than Theorem \ref{thm2}. In the following sections we make use of both Theorems.
\section{Multiplicativity of Renormalization Function with Complex Scaling}
\label{r3}
In this section we study the renormalization of summatory function when the 
summation variable is scaled by a product $\prod_{i=1}^np_i^{m_i}$ with $n$ 
distinct primes $p_i\geq 2$. Consider $F\left\{f;N,\prod_{i=1}^np_i^{m_i}
\right\}$ and find its governing functional equation (\ref{q66}).

First, write the relationship between two summatory functions with two different
summands, $f\left(k_1p_n^{m_n}\right)$ and $f(k_1)$, where $k_1=k\prod_{i=1}^
{n-1}p_i^{m_i}$. It is similar to that given in (\ref{b4}) and follows from the 
latter by replacing $k\rightarrow k_1$, i.e.,
\bea
F\left\{f;N,\prod_{i=1}^np_i^{m_i}\right\}=\sum_{r_n=0}^{p_n^{r_n}\leq N}L_{r_n}
\left(f;p_n^{m_n}\right)F\left\{f;\left\lfloor\frac{N}{p_n^{r_n}}\right\rfloor,
\prod_{i=1}^{n-1}p_i^{m_i}\right\}\;.\label{q30}
\eea
Next, repeat this procedure to reduce the scale by $p_{n-1}^{r_{n-1}}$ for 
summatory function appeared in the r.h.s. of (\ref{q30}) and substitute it again
into (\ref{q30}),
\bea
F\left\{f;N,\prod_{i=1}^np_i^{m_i}\right\}=\sum_{r_{n-1},r_n=0}^{p_{n-1}^{n-1}
p_n^{r_n}\leq N}L_{r_n}\left(f;p_n^{m_n}\right)L_{r_{n-1}}\left(f;p_{n-1}^{m_{
n-1}}\right)F\left\{f;\left\lfloor\frac{N}{p_{n-1}^{r_{n-1}}p_n^{r_n}}\right
\rfloor,\prod_{i=1}^{n-2}p_i^{m_i}\right\}\;.\nonumber
\eea
Continue to reduce the scales in a consecutive way for the next summatories and 
get finally,
\bea
F\left\{f;N,\prod_{i=1}^np_i^{m_i}\right\}=\sum_{r_1,\ldots,r_n=0}^{\prod_{i=1}
^np_i^{r_i}\leq N}\left(\prod_{i=1}L_{r_i}\left(f;p_i^{m_i}\right)\right)
F\left\{f;\left\lfloor\frac{N}{\prod_{i=1}^np_i^{r_i}}\right\rfloor,1\right\}\;.
\label{q66}
\eea
Define new renormalization functions,
\bea
R\left(f;N,\prod_{i=1}^np_i^{m_i}\right)=\frac{F\left\{f;N,\prod_{i=1}^np_i^{m_i
}\right\}}{F\{f;N,1\}}\;,\quad R_{\infty}\left(f;\prod_{i=1}^np_i^{m_i}\right):=
\lim_{N\to\infty}R\left(f;N,\prod_{i=1}^np_i^{m_i}\right)\quad\label{q65}
\eea
and find their representations through the characteristic functions $L_r\left(
f;p_i^{m_i}\right)$, $i=1,\ldots,n$, and degrees $a_1,b_1$ and $a_2,b_2$ of 
leading asymptotics $G_1(N)$ and $G_2(N)$, respectively. Following an approach 
developed in section \ref{r2}, represent $R\left(f;N,\prod_{i=1}^np_i^{m_i}
\right)$ as a sum
\bea
R\left(f;N,\prod_{i=1}^np_i^{m_i}\right)={\cal R}_1\left(f;N,\prod_{i=1}^np_i^{
m_i}\right)+{\cal R}_2\left(f;N,\prod_{i=1}^np_i^{m_i}\right)\;,\label{q67}
\eea
where ${\cal R}_j\left(f;N,\prod_{i=1}^np_i^{m_i}\right)$, $j=1,2$, are 
analogous to those given in (\ref{b9}) and (\ref{b9b}). Substitute there $G_1(
N)$ and $G_2(N)$ given in (\ref{i7}), and obtain formulas analogous to 
those given in (\ref{b10a}) and (\ref{b10w}), (\ref{b10x}). Here they are
\bea
{\cal R}_1\left(f;N,\prod_{i=1}^np_i^{m_i}\right)=\sum_{r_1,\ldots,r_n=0}^{ 
\prod_{i=1}^np_i^{r_i}\leq N}\prod_{i=1}^n\left(\frac{L_{r_i}\left(f;p_i^{m_i}
\right)}{p_i^{a_1r_i}}\right)\left(1-\frac1{\log_cN}\sum_{i=1}^nr_i\log_c p_i
\right)^{b_1},\label{q68}
\eea
where a base $c$ is choosen in such a way that $2\leq c<\min\{p_1,\ldots,p_n\}$,
so that the upper summation bound $\prod_{i=1}^np_i^{m_i}\leq N$ is 
correspondent to inequality, $\sum_{j=1}^nr_i\log_c p_i\leq\log_c N$. Regarding 
${\cal R}_2\left(f;N,\prod_{i=1}^np_i^{m_i}\right)$, we have an upper bound
\bea
\left|{\cal R}_2\left(f;N,\prod_{i=1}^np_i^{m_i}\right)\right|\leq\frac{{\cal C}
}{N^{a_2}\;\left(\log_c N\right)^{b_2}}\;{\cal R}_3\left(f;N,\prod_{i=1}^np_i^{
m_i}\right)\;,\quad\mbox{where}\label{q69}
\eea
\bea
{\cal R}_3\left(f;N,\prod_{i=1}^np_i^{m_i}\right)=\sum_{r_1,\ldots,r_n=0}^{  
\prod_{i=1}^np_i^{r_i}\leq N}\prod_{i=1}^n\left(\frac{\left|L_{r_i}\left(f;p_i^{
m_i}\right)\right|}{p_i^{r_i(a_1-a_2)}}\right)\left(1-\frac1{\log_cN}\sum_{i=1}
^nr_i\log_cp_i\right)^{b_1-b_2}.\quad\label{q70}
\eea
Keeping in mind the {\sf SQ} theorem and its usage in sections \ref{r21} and 
\ref{r22} we assume throughout this section $b_1,b_2\in{\mathbb Z}$. Extension 
on non integers $b_1$ and $b_2$ is trivial and can be done following those given
in Lemmas \ref{lem1} and \ref{lem2}, and therefore will be skipped. In next 
sections we prove several statements on ${\cal R}_1\left(f;N,\prod_{i=1}^np_i^{
m_i}\right)$ and ${\cal R}_2\left(f;N,\prod_{i=1}^np_i^{m_i}\right)$ which are 
similar to Lemmas \ref{lem1}, \ref{lem2}, \ref{lem3} and \ref{lem4} in section 
\ref{r2}. In this conjunction, it is important to use the same sufficient 
conditions which were used in these Lemmas.
\subsection{Convergence of ${\cal R}_1\left(f;N,\prod_{i=1}^np_i^{m_i}\right)$}
\label{r31}
\begin{lemma}\label{lem5}Let a function $f(k)\in{\mathbb B}\left\{N^{a_1}\left(
\log N\right)^{b_1};\;G_2(N)\right\}$, $a_1\geq 0$, $b_1\in{\mathbb Z}_+\cup
\{0\}$, be given and let there exist two numbers ${\cal K}>0$ and $\gamma<a_1$ 
and an integer $r_*\geq 0$ such that $|L_r\left(f;p^m\right)|\leq {\cal K}p^{
\gamma r}$ for all $r\geq r_*$. Then
\bea
\lim_{N\to\infty}{\cal R}_1\left(f;N,\prod_{i=1}^np_i^{m_i}\right)=\prod_{i=1}^
n\left(\sum_{r_i=0}^{\infty}\frac{L_{r_i}\left(f;p_i^{m_i}\right)}{p_i^{a_1r_i}}
\right)\;.\label{q71}
\eea
\end{lemma}
{\sf Proof} $\;\;\;$Exponentiating the binomial in (\ref{q68}) we obtain
\bea
{\cal R}_1\left(f;N,\prod_{i=1}^np_i^{m_i}\right)\!=\!\!\sum_{r_1,\ldots,r_n=0}
^{\prod_{i=1}^np_i^{r_i}\leq N}\!\!\left(\prod_{i=1}^n\frac{L_{r_i}\left(f;p_i^{
m_i}\right)}{p_i^{a_1r_i}}\right)\!+\sum_{k=1}^{b_1}(-1)^k{b_1\choose k}
{\cal R}_6\left(f;N,\prod_{i=1}^np_i^{m_i},k\right),\quad\label{q72}
\eea
where
\bea
{\cal R}_6\left(f;N,\prod_{i=1}^np_i^{m_i},k\right)=\frac1{\left(\log_cN\right)
^k}\sum_{r_1,\ldots,r_n=0}^{\prod_{i=1}^np_i^{r_i}\leq N}\left(\prod_{i=1}^n
\frac{L_{r_i}\left(f;p_i^{m_i}\right)}{p_i^{a_1r_i}}\right)\;\left(\sum_{i=1}^n
r_i\log_cp_i\right)^k.\nonumber
\eea
Find an estimate for ${\cal R}_6\left(f;N,\prod_{i=1}^np_i^{m_i},k\right)$,
\bea
\left|{\cal R}_6\left(f;N,\prod_{i=1}^np_i^{m_i},k\right)\right|\leq
{\cal K}^n\!\!\sum_{k_1,\ldots,k_n=0\atop k_1+\ldots+k_n=k}^k\!{k\choose k_1,
\ldots,k_n}\prod_{i=1}^n\left(\frac{\log_cp_i}{\log_cN}\right)^{k_i}\;
\sum_{r_1,\ldots,r_n=0}^{\prod_{i=1}^np_i^{r_i}\leq N}
\left(\prod_{i=1}^n\frac{r_i^{k_i}}{p_i^{\epsilon r_i}}\right),\nonumber
\eea
where $\epsilon=a_1-\gamma>0$. One more inequality reads
\bea
\sum_{r_1,\ldots,r_n=0}^{\prod_{i=1}^np_i^{r_i}\leq N}\left(\prod_{i=1}^n\frac{
r_i^{k_i}}{p_i^{\epsilon r_i}}\right)\leq\prod_{i=1}^n\left(\sum_{r_i=0}^{p_i^{
r_i}\leq N}\frac{r_i^{k_i}}{p_i^{\epsilon r_i}}\right)=\prod_{i=1}^nT(p_i,k_i,
\epsilon,\log_{p_i}N)\;,\label{q42}
\eea
where $T(p,k,\epsilon,M)$ is defined in (\ref{b11c1}). Combining the two last
inequalities together we get
\bea
\left|{\cal R}_6\left(f;N,\prod_{i=1}^np_i^{m_i},k\right)\right|\leq{\cal K}^n
\sum_{k_1,\ldots,k_n=0\atop k_1+\ldots+k_n=k}^k\!{k\choose k_1,\ldots,k_n}
\prod_{i=1}^n\frac{T(p_i,k_i,\epsilon,\log_{p_i}N)}{\left(\log_{p_i}N\right)^{
k_i}}\;.\label{q73}
\eea
Inserting the asymptotics (\ref{b11b1}) of $T(p_i,k_i,\epsilon,\log_{p_i}N)$
into (\ref{q73}) we arrive at
\bea
\left|{\cal R}_6\left(f;N,\prod_{i=1}^np_i^{m_i},k\right)\right|\leq{\cal K}^n
\sum_{k_1,\ldots,k_n=0\atop k_1+\ldots+k_n=k}^k\!{k\choose k_1,\ldots,k_n}
\prod_{i=1}^n\left[\frac{{\rm Li}_{-k_i}\left(p_i^{-\epsilon}\right)}{\left(
\log_{p_i}N\right)^{k_i}}- \frac{N^{-\epsilon}}{p_i^{\epsilon}-1}\right]\;.
\label{q74}
\eea
Repeating the concluding remarks in proof of Lemma \ref{lem1} on asymptotics of
the polylogarithm function ${\rm Li}_s(z)$ we conclude that ${\cal R}_6\left(f;N
,\prod_{i=1}^np_i^{m_i},k\right)$ is convergent to zero when $N\to\infty$. Then,
keeping in mind the representation (\ref{q72}) for ${\cal R}_1\left(f;N,\prod_{
i=1}^np_i^{m_i}\right)$ and running the upper bound of summation to infinity we 
conclude that the limit (\ref{q71}) holds.$\;\;\;\;\;\;\Box$
\begin{lemma}\label{lem6}Let a function $f(k)\in{\mathbb B}\left\{N^{a_1}\left(
\log N\right)^{b_1};\;G_2(N)\right\}$, $a_1>0$, $b_1\in{\mathbb Z}_-$, be given 
and let there exist two numbers ${\cal K}>0$ and $\gamma<a_1$ and an integer 
$r_*\geq 0$ such that $|L_r\left(f;p^m\right)|\leq {\cal K}p^{\gamma r}$ for all
$r\geq r_*$. Then (\ref{q71}) holds.
\end{lemma}
{\sf Proof} $\;\;\;$Here we follow the keyline in the proof of Lemma \ref{lem2} 
and, according to (\ref{b6c}) in Remark \ref{rem2}, start with representation of
${\cal R}_1\left(f;N,\prod_{i=1}^np_i^{m_i}\right)$ avoiding its divergence at 
$\prod_{i=1}^np_i^{m_i}=N$,
\bea
{\cal R}_1\left(f;N,\prod_{i=1}^np_i^{m_i}\right)=\sum_{r_1,\ldots,r_n=0}^{
\prod_{i=1}^np_i^{m_i}\leq N-1}\left(\prod_{i=1}^n\frac{L_{r_i}\left(f;p_i^{m_i}
\right)}{p_i^{a_1r_i}}\right)\left(1-\frac1{\log_cN}\sum_{i=1}^nr_i\log_cp_i
\right)^{-|b_1|}.\quad\label{q44}
\eea
Making use of identity (\ref{b218a}) we obtain,
\bea
{\cal R}_1\left(f;N,\prod_{i=1}^np_i^{m_i}\right)=\sum_{r_1,\ldots,
r_n=0}^{\prod_{i=1}^np_i^{m_i}\leq N-1}\!\!\left(\prod_{i=1}^n\frac{L_{r_i}
\left(f;p_i^{m_i}\right)}{p_i^{a_1r_i}}\right)+\sum_{k=1}^{|b_1|}{\cal R}_7
\left(f;N,\prod_{i=1}^np_i^{m_i},k\right),\quad\label{q44a}
\eea
where $\tilde{r}=\sum_{i=1}^nr_i\log_cp_i$ and 
\bea
{\cal R}_7\left(f;N,\prod_{i=1}^np_i^{m_i},k\right)=\frac1{\log_cN}\sum_{r_1,
\ldots,r_n=0}^{\prod_{i=1}^np_i^{m_i}\leq N-1}\left(\prod_{i=1}^n\frac{L_{r_i}
\left(f;p_i^{m_i}\right)}{p_i^{a_1r_i}}\right)\;\tilde{r}\left(1-
\frac{\tilde{r}}{\log_cN}\right)^{-k}.\label{q44b}
\eea
Making use of inequality (\ref{b20i}) and constraint imposed on $L_{r}\left
(f;p^{m}\right)$
\bea
\left(1-\frac{\tilde{r}}{\log_cN}\right)^{-k}\leq (1+\tilde{r})^k\;,\quad
|L_r\left(f;p^m\right)|\leq {\cal K}p^{\gamma r}\;,\nonumber
\eea
exponentiate the binomial $(1+\tilde{r})^k$ in (\ref{q44b}) and obtain
\bea
\left|{\cal R}_7\left(f;N,\prod_{i=1}^np_i^{m_i},k\right)\right|\leq\frac{{\cal 
K}^n}{\log_c N}\!\sum_{r_1,\ldots,r_n=0}^{\prod_{i=1}^np_i^{m_i}\leq N-1}\!
\frac{\tilde{r}(1+\tilde{r})^k}{\prod_{i=1}^np_i^{\epsilon r_i}}=\frac{{\cal 
K}^n}{\log_c N}\sum_{j=0}^k{k\choose j}\sum_{r_1,\ldots,r_n=0}^{\prod_{i=1}^n
p_i^{m_i}\leq N-1}\!\frac{\tilde{r}^{j+1}}{p_1^{\epsilon r_1}p_2^{\epsilon r_2}}
.\nonumber
\eea
where $\epsilon=a_1-\gamma$. Exponentiating a binomial $\left(\sum_{i=1}^nr_i
\log_c p_i\right)^{j+1}$ in the last expression we get
\bea
\left|{\cal R}_7\left(f;N,p_1^{m_1}p_2^{m_2},k\right)\right|\leq\frac{{\cal K}^n
}{\log_c N}\sum_{j=0}^k{k\choose j}\sum_{j_1,\ldots,j_n=0\atop j_1+\ldots+j_n=
j+1}^{j+1}{j+1\choose j_1,\ldots,j_n}\sum_{r_1,\ldots,r_n=0}^{\prod_{i=1}^np_i^{
m_i}\leq N-1}\prod_{i=1}^n\frac{r_i^{j_i}}{p_i^{\epsilon r_i}}\left(\log_c p_i
\right)^{j_i},\nonumber
\eea
Asymptotic behavior in $N$ of the last expression is completely determined by 
its inner sum in $r_i$ with respect to its prefactor $\left(\log_c N\right)^{
-1}$. This behavior can be calculated following corresponding part (\ref{q42}) 
of the proof in Lemma \ref{lem1}, 
\bea
\sum_{r_1,\ldots,r_n=0}^{\prod_{i=1}^np_i^{m_i}\leq N-1}\!\prod_{i=1}
\frac{r_i^{j_i}}{p_i^{\epsilon r_i}}\left(\log_c p_i\right)^{j_i}\stackrel{
N\to\infty}{\leq}\prod_{i=1}\left[{\rm Li}_{-j_i}\left(p_i^{-\epsilon}\right)-
\frac{N^{-\epsilon}\left(\log_{p_i}N\right)^{j_i}}{1-p_i^{-\epsilon}}\right]
\left(\log_c p_i\right)^{j_i}.\quad\label{q46}
\eea
Keeping in mind the prefactor $\left(\log_cN\right)^{-1}$ and the last 
asymptotics (\ref{q46}), the upper bound for ${\cal R}_7\left(f;N,\prod_{i=1}^
np_i^{m_i},k\right)$ can be done infinitely small, i.e., it is convergent to 
zero when $N\to\infty$.

Making use of representation (\ref{q44a}) for ${\cal R}_1\left(f;N,p_1^{m_1}p_2
^{m_2}\right)$ and running the upper bound of summation to infinity we conclude
that the limit (\ref{q71}) holds.$\;\;\;\;\;\;\Box$
\vspace{-.3cm}
\subsection{Convergence of ${\cal R}_2\left(f;N,\prod_{i=1}^np_i^{m_i}\right)$}
\label{r32}
\begin{lemma}\label{lem7}Let a function $f(k)\in{\mathbb B}\left\{N^{a_1}\left(
\log N\right)^{b_1};N^{-a_2}\left(\log N\right)^{-b_2}\right\}$, $b_1\geq b_2$,
be given and let there exist two numbers ${\cal K}>0$ and $\gamma<a_1$ and an
integer $r_*\geq 0$ such that $|L_r\left(f;p^m\right)|\leq {\cal K}p^{\gamma r}$
for all $r\geq r_*$, then
\bea
\lim_{N\to\infty}{\cal R}_2\left(f;N,\prod_{i=1}^np_i^{m_i}\right)=0\;.
\label{q75}
\eea
\end{lemma}
{\sf Proof} $\;\;\;$Denote $b_1-b_2=e\in{\mathbb Z}_+\cup\{0\}$ and make use of
an inequality, $(1-x)^e\leq 1$ when $0\leq x\leq 1$, $e\geq 0$. Substituting    
constraints on $L_{r}\left(f;p^{m}\right)$ into (\ref{q69}) we obtain,
\bea
\left|{\cal R}_2\left(f;N,\prod_{i=1}^np_i^{m_i}\right)\right|\leq\frac{{\cal C}
{\cal K}^n\;{\cal I}_n(\nu)}{N^{a_2}\left(\log_p N\right)^{b_2}},\quad\mbox{
where}\quad{\cal I}_n(\nu)=\sum_{r_1,\ldots,r_n=0}^{\prod_{i=1}^np_i^{r_i}\leq
N}\prod_{i=1}^np_i^{\nu r_i}\;.\label{q76}
\eea
and $\nu=a_2-a_1+\gamma$. Focus on the sum ${\cal I}_n(\nu)$ and estimate it in 
3 cases, $\nu=0$, $\nu>0$ and $\nu<0$. 

Let $\nu=0$, i.e., $a_2=a_1-\gamma$, then ${\cal I}_n(0)$ accounts for a number
of integral points (vertices with integer coordinates) in the $n$-dim simplex, 
or corner of the $n$-dim cube, defined as follows,
\bea
\Delta_n=:\left\{r_1,\ldots,r_n\in{\mathbb Z}_+\cup\{0\}\;|\;\sum_{i=1}^nr_i
\log_cp_i\leq \log_cN\right\}\;.\nonumber
\eea
The simplex $\Delta_n$ has one orthogonal corner and sizes of edges $\log_cN/
\log_cp_i$ along the $i$th axis. The number ${\cal I}_n(0)$ is described by the
Ehrhart polynomial \cite{dr97} and, when $N\to\infty$, it has a leading term 
coinciding with simplex' volume.
\bea
{\cal I}_n(0)=\sum_{r_1,\ldots,r_n=0}^{\prod_{i=1}^np_i^{r_i}\leq N}\simeq 
\frac{\left(\log_c N\right)^n}{n!\;\prod_{i=1}\log_cp_i}\label{q77}
\eea
Substituting (\ref{q77}) into (\ref{q76}) we get
\bea
\left|{\cal R}_2\left(f;N,\prod_{i=1}^np_i^{r_i}\right)\right|\leq
\frac{{\cal C}{\cal K}^n}{n!\;\prod_{i=1}\log_cp_i}\;
\frac{N^{\gamma-a_1}}{\left(\log_c N\right)^{b_2-n}}\;.\label{q50}
\eea
By (\ref{q50}) and inequality $\gamma<a_1$ we conclude that ${\cal R}_2\left(f;
N,p_1^{m_1}p_2^{m_2}\right)$ is convergent to zero when $N\to\infty$ 
irrespectively to the value of $b_2$. 

Consider the case $\nu>0$ and estimate ${\cal I}_n(\nu)$ and ${\cal R}_2\left(
f;N,\prod_{i=1}^np_i^{r_i}\right)$,
\bea
{\cal I}_n(\nu)\leq N^{\nu}{\cal I}_n(0)\;,\quad\rightarrow\quad
\left|{\cal R}_2\left(f;N,\prod_{i=1}^np_i^{r_i}\right)\right|\leq
\frac{{\cal C}{\cal K}^n}{n!\;\prod_{i=1}\log_cp_i}\;
\frac{N^{\gamma-a_1}}{\left(\log_c N\right)^{b_2-n}}\;.\label{q51}
\eea
By (\ref{q51}) the term ${\cal R}_2\left(f;N,\prod_{i=1}^np_i^{r_i}\right)$ is 
also convergent to zero when $N\to\infty$. 

Finally, consider the case $\nu<0$. According to (\ref{q76}) we have 
${\cal I}_n(\nu)<{\cal I}_n(0)$ and therefore
\bea
\left|{\cal R}_2\left(f;N,\prod_{i=1}^np_i^{r_i}\right)\right|\leq
\frac{{\cal C}{\cal K}^n}{n!\;\prod_{i=1}\log_cp_i}\;
\frac{N^{-a_2}}{\left(\log_c N\right)^{b_2-n}}\;.\label{q52}
\eea
Thus, ${\cal R}_2\left(f;N,\prod_{i=1}^np_i^{r_i}\right)$ is convergent to zero 
due to (\ref{q52}) and constraints (\ref{i7}) on degrees $a_2$ and $b_2$. 
Summarizing (\ref{q50}), (\ref{q51}) and (\ref{q52}) we complete the proof 
of Lemma.$\;\;\;\;\;\;\Box$

In the following Lemmas we consider two different cases, $a_2>0$ and $a_2=0$, 
separately.
\begin{lemma}\label{lem8}Let a function $f(k)\in{\mathbb B}\left\{N^{a_1}\left(
\log N\right)^{b_1};N^{-a_2}\left(\log N\right)^{-b_2}\right\}$, $b_1<b_2$ and 
$a_2>0$, be given and let there exist two numbers ${\cal K}>0$ and $\gamma<a_1$ 
and an integer $r_*\geq 0$ such that $|L_r\left(f;p^m\right)|\leq {\cal K}p^{
\gamma r}$ for all $r\geq r_*$. Then (\ref{q75}) holds.
\end{lemma}
{\sf Proof} $\;$In accordance with (\ref{b6c}) in Remark \ref{rem2}, represent 
${\cal R}_3\left(f;N,\prod_{i=1}^np_i^{m_i}\right)$ avoiding its divergence at 
$\prod_{i=1}^np_i^{m_i}=N$,
\bea
{\cal R}_3\left(f;N,\prod_{i=1}^np_i^{m_i}\right)=\sum_{r_1,\ldots,r_n=0}^{
\prod_{i=1}^np_i^{m_i}\leq N-1}\left(\prod_{i=1}^n\frac{L_{r_i}\left(f;p_i^{ 
m_i}\right)}{p_i^{r_i(a_1-a_2)}}\right)\left(1-\frac{\tilde{r}}{\log_c N}
\right)^{-(b_2-b_1)},\label{q55a}
\eea
and make use of inequality (\ref{b20i}) in the range $0\leq\tilde{r}\leq\left
\lfloor\log_p N\right\rfloor-1$,
\bea
\left(1-\frac{\tilde{r}}{\log_c N}\right)^{-(b_2-b_1)}\leq\left(\log_cN\right)
^{b_2-b_1}\;.\label{q55}
\eea
Substituting (\ref{q55}) into (\ref{q55a}) for ${\cal R}_3\left(f;N,\prod_{i=1}
^np_i^{m_i}\right)$ we get
\bea
{\cal R}_3\left(f;N,\prod_{i=1}^np_i^{m_i}\right)\leq\left(\log_c N\right)^{b_2-
b_1}\sum_{r_1,\ldots,r_n=0}^{\prod_{i=1}^np_i^{m_i}\leq N-1}\prod_{i=1}^n\frac{
\left|L_{r_i}\left(f;p_i^{m_i}\right)\right|}{p_i^{r_1(a_1-a_2)}}\;.\label{q56}
\eea
Apply to (\ref{q56}) the constraints on $L_{r}\left(f;p^{m}\right)$ and 
substitute the result into (\ref{q69}),
\bea
{\cal R}_2\left(f;N,\prod_{i=1}^np_i^{m_i}\right)\leq\frac{{\cal C}{\cal K}^n}
{N^{a_2}\;\left(\log_cN\right)^{b_1}}\sum_{r_1,\ldots,r_n=0}^{\prod_{i=1}^np_i^{
m_i}\leq N-1}\prod_{i=1}^np_i^{\nu r_i}\;,\quad \nu=a_2-a_1+\gamma\;.\label{q57}
\eea
By comparison (\ref{q57}) with (\ref{q76}) from Lemma \ref{lem7} we obtain
according to (\ref{q50}), (\ref{q51}) and (\ref{q52})
\bea
\left|{\cal R}_2\left(f;N,\prod_{i=1}^np_i^{m_i}\right)\right|&\leq&
\frac{{\cal C}{\cal K}^n}{n!\;\prod_{i=1}\log_cp_i}\;\frac{N^{\gamma-a_1}}
{\left(\log_c N\right)^{b_1-n}}\;,\quad \nu\geq 0\;,\label{q58}\\
\left|{\cal R}_2\left(f;N,\prod_{i=1}^np_i^{m_i}\right)\right|&\leq&
\frac{{\cal C}{\cal K}^n}{n!\;\prod_{i=1}\log_cp_i}\;\frac{N^{-a_2}}
{\left(\log_c N\right)^{b_1-n}}\;,\quad \nu<0\;.\label{q59}
\eea
Thus, by (\ref{q58}) and inequality $\gamma<a_1$ a series ${\cal R}_2\left(f;N,
\prod_{i=1}^np_i^{m_i}\right)$ is convergent to zero when $N\to\infty$ 
irrespectively to the value of $b_1$. The same conclusion (\ref{q75}) on 
convergence of this series holds due to (\ref{q59}) when $a_2>0$.
$\;\;\;\;\;\;\Box$
\begin{lemma}\label{lem9}Let a function $f(k)\in{\mathbb B}\left\{N^{a_1}\left(
\log N\right)^{b_1};\left(\log N\right)^{-b_2}\right\}$, $b_1<b_2$, $a_2=0$, be 
given and let there exist two numbers ${\cal K}>0$ and $\gamma<a_1$ and an 
integer $r_*\geq 0$ such that $|L_r\left(f;p^m\right)|\leq {\cal K}p^{\gamma r}$
for all $r\geq r_*$. Then (\ref{q75}) holds.
\end{lemma}
{\sf Proof} $\;\;\;$This case has to be treated more precisely than that in 
Lemma \ref{lem8}. Rewrite (\ref{q55a})
\bea
{\cal R}_3\left(f;N,\prod_{i=1}^np_i^{m_i}\right)=\sum_{r_1,\ldots,r_n=0}^{
\prod_{i=1}^np_i^{m_i}\leq N-1}\prod_{i=1}^n\frac{\left|L_{r_i}\left(f;p_i^{
m_i}\right)\right|}{p_i^{a_1r_i}}\left(1-\frac{\tilde{r}}{\log_cN}\right)^{-(
b_2-b_1)},\label{q61}
\eea
and compare it with expression (\ref{q44}) when $b_1<0$. A difference in 
degrees, $b_2-b_1$ and $b_1$, does not break the main result of Lemma \ref{lem6}
: only the 1st leading term in (\ref{q44a}) is survived when $N\to\infty$. When
we apply it to (\ref{q61}) and make use of constraint on $L_{r}\left(f;p^{m}
\right)$ we get
\bea
{\cal R}_3\left(f;N,\prod_{i=1}^np_i^{m_i}\right)\stackrel{N\to\infty}{=}\sum_{
r_1,\ldots,r_n=0}^{\prod_{i=1}^np_i^{m_i}\leq N-1}\prod_{i=1}^n\frac{L_{r_i}
\left(f;p_i^{m_i}\right)}{p_i^{a_1r_i}}\leq{\cal K}^n\prod_{i=1}^n\left(\sum_{
r_i=0}^{\infty}p_i^{-\epsilon r_i}\right)=\frac{{\cal K}^n}{\prod_{i=1}^n\left(
1-p_i^{-\epsilon}\right)},\nonumber
\eea
where $\epsilon=a_1-\gamma>0$. Substituting the last estimate into (\ref{q69})
we obtain,
\bea
\left|{\cal R}_2\left(f;N,\prod_{i=1}^np_i^{m_i}\right)\right|\leq\frac{
{\cal C}{\cal K}^n}{\prod_{i=1}^n\left(1-p_i^{-\epsilon}\right)}\cdot\frac1{
\left(\log_p N\right)^{b_2}}\;.\label{q62}
\eea
Recall that by (\ref{i7}) the degrees of the error term satisfy the condition:
if $a_2=0$ than $b_2>0$. Thus, by (\ref{q62}) the series ${\cal R}_2\left(f;N,
\prod_{i=1}^np_i^{m_i}\right)$ converges to zero when $N\to\infty$ that proves 
Lemma.$\;\;\;\;\;\;\Box$

We combine Lemmas \ref{lem5}, \ref{lem6}, \ref{lem7}, \ref{lem8} and \ref{lem9}
on convergence of ${\cal R}_j\left(f;N,\prod_{i=1}^np_i^{m_i}\right)$, $j=1,2$,
and according to (\ref{q67}) we arrive at the analogue of Theorem \ref{thm1} in
the case of a scaling by $\prod_{i=1}^np_i^{m_i}$.
\begin{theorem}\label{thm3}Let a function $f(k)\in{\mathbb B}\left\{G_1(N);G_2
(N)\right\}$ be given with $a_i$ and $b_i$ satisfying (\ref{i7}) and let there
exist two numbers ${\cal K}>0$ and $\gamma<a_1$ and an integer $r_*\geq 0$ such
that $|L_r\left(f;p^m\right)|\leq {\cal K}p^{\gamma r}$ for all $r\geq r_*$.
Then
\bea
R_{\infty}\left(f;\prod_{i=1}^np_i^{m_i}\right)=\prod_{i=1}^n\left(\sum_{r_i=0}
^{\infty}\frac{L_{r_i}\left(f;p_i^{m_i}\right)}{p_i^{a_1r_i}}\right)\;.
\label{q63}
\eea
\end{theorem}
Combining Theorems \ref{thm1} and \ref{thm3} we come to important consequence 
which manifests the multiplicative property of the renormalization function.
\begin{corollary}\label{cor1}
Under the conditions of Theorem \ref{thm3} the following holds,
\bea
R_{\infty}\left(f;\prod_{i=1}^np_i^{m_i}\right)=\prod_{i=1}^nR_{\infty}\left(f;
p_i^{m_i}\right)\;.\label{q64}
\eea
\end{corollary}
\subsection{Asymptotics of Summatory Functions $\sum_{k_1,k_2\leq N}f(k_1k_2)$}
\label{r33}
In this section we calculate the summatory function $\Phi[f;N,1]=\sum_{k_1,k_2
\leq N}f(k_1k_2)$ and find its asymptotics by applying Corollary \ref{cor1}. 
According to definition of summatory function we get,
\bea
\Phi\{f;N,1\}=\!\sum_{k\leq N}f(k)+\sum_{k\leq N}f(2k)+\ldots=\sum_{k\leq N}
\sum_{m_{i,j}=0\atop n_j=1}\!f\left(k\prod_{i=1}^{n_j}p_i^{m_{i,j}}\right)\!=\!
\sum_{m_{i,j}=0\atop n_j=1}\sum_{k\leq N}f\left(k\prod_{i=1}^{n_j}p_i^{m_{i,j}}
\right)\nonumber
\eea
where indices $i$, $j$, $m_{i,j}$ and $n_j$ account for all primes such that 
$\prod_{i=1}^{n_j}p_i^{m_{i,j}}\leq N$. Thus, according to definition of 
summatory function with scaled summation variable, we obtain,
\bea
\Phi\{f;N,1\}=\sum_{m_{i,j}=0\atop n_j=1}F\left\{f;N,\prod_{i=1}^{n_j}p_i^{m_{
i,j}}\right\},\quad\frac{\Phi\{f;N,1\}}{F\{f;N,1\}}=\sum_{m_{i,j}=0\atop n_j=1}
R\left(f;N,\prod_{i=1}^{n_j}p_i^{m_{i,j}}\right),\label{n1} 
\eea
where $\prod_{i=1}^{n_j}p_i^{m_{i,j}}\leq N$. Consider asymptotics (omitting the
error terms) of three summatory functions when $N\to\infty$,
\bea
&&F\{f;N,1\}\stackrel{N\to\infty}{=}{\cal F}G_1(N)\;,\hspace{1.5cm}\Phi\{f;N,1
\}\stackrel{N\to\infty}{=}{\cal G}_1(f)\Gamma_1(N),\nonumber\\
&&\sum_{m_{i,j}=0\atop n_j=1}R\left(f;N,\prod_{i=1}^{n_j}p_i^{m_{i,j}}\right)
\stackrel{N\to\infty}{=}{\cal G}_2(f)\Gamma_2(N)\;.\label{n2}
\eea
Combining (\ref{n2}) and the 2nd formula in (\ref{n1}) we obtain
\bea
\Gamma_1(N)=G_1(N)\cdot\Gamma_2(N)\;,\quad{\cal G}_1(f)={\cal F}\cdot
{\cal G}_2(f)\;.\label{n3}
\eea
Calculation of $\Gamma_2(N)$ and ${\cal G}_2(f)$ is a difficult numerical task. 
Consider a special case when ${\cal G}_2(f)$ may be given in a closed form, 
namely, when $\Gamma_2(N)=N^0$, i.e., $\Gamma_1(N)=G_1(N)$. Consider the 3rd 
asymptotics in (\ref{n2}) and, according to Corollary \ref{cor1}, find its 
limit when $N\to\infty$,
\bea
{\cal G}_2(f)&=&\lim_{N\to\infty}\sum_{m_{i,j}=0\atop n_j=1}^{\prod_{i=1}^{n_j}
p_i^{m_{i,j}}\leq N}R\left(f;N,\prod_{i=1}^{n_j}p_i^{m_{i,j}}\right)=
\lim_{N\to\infty\atop n_j\to\infty}\sum_{m_{i,j}=0\atop n_j=1}^{\log_{p_i} N}
\prod_{i=1}^{n_j}R_{\infty}\left(f;p_i^{m_{i,j}}\right)\nonumber\\
&=&\lim_{n\to\infty}\prod_{i=1}^n\sum_{m=0}^{\infty}R_{\infty}\left(f;p_i^{m}
\right)=\prod_{p\geq 2}\;\sum_{m=0}^{\infty}R_{\infty}\left(f;p^m\right)\;,
\label{n5}
\eea
Another special case comes when $f(k)$ is a completely multiplicative arithmetic
function. i.e., $f(k_1k_2)=f(k_1)f(k_2)$. This leads to equalities: $\Gamma_2(N)
\equiv\Gamma_1(N)$ and ${\cal G}_2(f)={\cal F}$. We will illustrate this 
statement and (\ref{n5}) in section \ref{r42}.
\section{Renormalization of Dirichlet Series and Others Summatory Functions}
\label{r4}
In this section we extend the renormalization approach on summatory functions 
of more complex structure. They involve the summatory functions with summands 
given by $\prod_{i=1}^nf_i(k)$ and summation variable $k$ scaled for every 
multiplicative function $f_i$ by $p^{m_i}$, $m_i\neq m_j$. The case of the
Dirichlet series is a special one when $n=2$ and $f_2(k)=k^{-s}$. We study also 
the renormalization of summatory functions with summands given by $f(k^n)$. 
\subsection{Renormalization of Summatory Function $\sum_{k\leq N}\prod_{i=1}^n
f_i\left(kp^{m_i}\right)$}\label{r41}
Start with summatory function $F\left\{\prod_{i=1}^nf_i;N,p^{{\bf m}}\right\}=
\sum_{k\leq N}\prod_{i=1}^nf_i\left(kp^{m_i}\right)$, where $p^{{\bf m}}$
denotes a tuple $\left\{p^{m_1},\ldots,p^{m_n}\right\}$, and make use of a 
standard notation $F\left\{\prod_{i=1}^nf_i;N,1^{{\bf m}}\right\}=F\left\{\prod
_{i=1}^nf_i;N,1\right\}$. Derive for $F$ a functional equation following the 
approach developed in section \ref{r12} and start
\bea
F\left\{\prod_{i=1}^nf_i;N,p^{{\bf m}}\right\}&=&{\sf F}\left(p^{{\bf m}}\right)
\sum_{k=1\atop p\;\nmid\;k}^N\prod_{i=1}^nf_i(k)+\sum_{k=p\atop 
p\;\mid\;k}^N\prod_{i=1}^nf_i\left(kp^{m_i}\right)\nonumber\\
&=&{\sf F}\left(p^{{\bf m}}\right)\left(\sum_{k=1}^N\prod_{i=1}^nf_i(k)-\sum_{
k=p\atop p\;\mid\;k}^N\prod_{i=1}^nf_i(k)\right)+\sum_{l=1}^{N_1}\prod_{i=1}^n
f_i\left(lp^{m_i+1}\right)\;,\nonumber
\eea
where 
${\sf F}\left(p^{{\bf m}}\right)=\prod_{i=1}^nf_i\left(p^{m_i}\right)$ and 
$N_r$ was defined in section \ref{r12}. Rewrite the last equality
\bea
F\left\{\prod_{i=1}^nf_i;N,p^{{\bf m}}\right\}-F\left\{\prod_{i=1}^nf_i;N_1,
p^{{\bf m+1}}\right\}={\sf F}\left(p^{{\bf m}}\right)\left[F\left\{\prod_{
i=1}^nf_i;N,1\right\}-F\left\{\prod_{i=1}^nf_i;N_1,p\right\}\right],\nonumber
\eea
which is similar to (\ref{b2}). The corresponding counterpartner for its general
version (\ref{b3}) reads,
\bea
F\left\{\prod_{i=1}^nf_i;N_r,p^{{\bf m}}\right\}\!-F\left\{\prod_{i=1}^nf_i;N_{
r+1},p^{{\bf m+1}}\right\}\!=\!{\sf F}\left(p^{{\bf m}}\right)\left[F\left\{
\prod_{i=1}^nf_i;N_r,1\right\}\!-F\left\{\prod_{i=1}^nf_i;N_{r+1},p\right\}
\right]\nonumber
\eea
Combining last equations of running index $0\leq r\leq\left\lfloor\log_p N\right
\rfloor$ together we arrive at the functional equation for summatory function,
\bea
F\left\{\prod_{i=1}^nf_i;N,p^{{\bf m}}\right\}=\sum_{r=0}^{\left\lfloor\log_p 
N\right\rfloor}L_r\left(\prod_{i=1}^nf_i;p^{{\bf m}}\right)F\left\{\prod_{i=1}^n
f_i;N_r,1\right\}\;,\label{n7}
\eea
where
\bea
L_r\left(\prod_{i=1}^nf_i;p^{{\bf m}}\right)=\prod_{i=1}^nf_i\left(p^{m_i+r}
\right)-\sum_{j=0}^{r-1}L_j\left(\prod_{i=1}^nf_i;p^{{\bf m}}\right)\prod_{i=1}^
nf_i\left(p^{r-j}\right)\;,\label{n8}
\eea
Formulas for the first $L_r\left(\prod_{i=1}^nf_i;p^{{\bf m}}\right)$ read
\bea
L_0\left(\prod_{i=1}^nf_i;p^{{\bf m}}\right)&=&\prod_{i=1}^nf_i\left(p^{m_i}
\right)\;,\quad 
L_1\left(\prod_{i=1}^nf_i;p^{{\bf m}}\right)=\prod_{i=1}^nf_i\left(p^{m_i+1}
\right)-\prod_{i=1}^nf_i\left(p^{m_i}\right)\prod_{i=1}^nf_i(p),\nonumber\\
L_2\left(\prod_{i=1}^nf_i;p^{{\bf m}}\right)&=&\prod_{i=1}^nf_i\left(p^{m_i+2}
\right)-\prod_{i=1}^nf_i\left(p^{m_i}\right)\prod_{i=1}^nf_i\left(p^2\right)-
\prod_{i=1}^nf_i\left(p^{m_i+1}\right)\prod_{i=1}^nf_i(p)+\nonumber\\
&&\prod_{i=1}^nf_i\left(p^{m_i}\right)\prod_{i=1}^nf_i^2(p)\;,\quad\mbox{etc}
\label{n9}
\eea
such that for $p=1$ or $m=0$ we have, $L_0\left(\prod_{i=1}^nf_i;1^{{\bf m}}
\right)=1$ and $L_r\left(\prod_{i=1}^nf_i;1^{{\bf m}}\right)=0$, $r\geq 1$. Find
analogues to formulas (\ref{b6a}) when $f_i(p^m)=A_{f_i}p^{m-1}$, $m\geq 1$, and
$A_{f_i}$ denotes the real constant. By (\ref{n8}) or (\ref{n9}) such formulas 
for $L_r\left(\prod_{i=1}^nf_i;p^{{\bf m}}\right)$ can be calculated by 
induction,
\bea
L_r\left(\prod_{i=1}^nf_i;p^{{\bf m}}\right)&=&p^{m_1+\ldots+m_n-n}
\left(p^n-\prod_{i=1}^nA_{f_i}\right)^r\prod_{i=1}^nA_{f_i}\;,\nonumber\\
L_r\left(\prod_{i=1}^n\frac1{f_i};p^{{\bf m}}\right)&=&p^{n-m_1-\ldots-m_n}
\left(p^{-n}-\prod_{i=1}^nA_{f_i}^{-1}\right)^r\prod_{i=1}^nA_{f_i}^{-1}
\;.\nonumber
\eea
E.g., in the case of the Euler $f_1=\phi(k)$ and Dedekind $f_2=\psi(k)$ totient 
functions we have
\bea
L_r\left(\phi\psi;p^{{\bf m}}\right)=(p^2-1)p^{m_1+m_2-2}\;,\quad 
L_r\left(\frac1{\phi\psi};p^{{\bf m}}\right)=\frac{(-1)^rp^{2(1-r)-m_1-m_2}}
{(p^2-1)^{r+1}}\;.\nonumber
\eea
Define new renormalization functions,
\bea
R\left(\prod_{i=1}^nf_i;N,p^{{\bf m}}\right)=\frac{F\left\{\prod_{i=1}^nf_i;N,
p^{{\bf m}}\right\}}{F\{\prod_{i=1}^nf_i;N,1\}}\;,\quad R_{\infty}\left(\prod_{
i=1}^nf_i;p^{{\bf m}}\right)=\lim_{N\to\infty}R\left(
\prod_{i=1}^nf_i;N,p^{{\bf m}}\right)\;.\label{n10}
\eea
By comparison formulas (\ref{n7}), (\ref{n8}), (\ref{n9}) with (\ref{b4}), 
(\ref{b5}), (\ref{b6}), respectively, and definition (\ref{n10}) with (\ref{i2})
we conclude that all results on renormalization of summatory function $F\left\{f
;N,p^m\right\}$ in section \ref{r2} can be reproduced for summatory function $F
\left\{\prod_{i=1}^nf_i;N,1^{{\bf m}}\right\}$ with a few necessary alterations.
Below, in Theorem \ref{thm4} we give (without proof) a sufficient condition for 
convergence of the asymptotics of renormalization function $R_{\infty}\left(
\prod_{i=1}^nf_i;p^{{\bf m}}\right)$. Its proof does not use new ideas and can 
be given following Lemmas \ref{lem1}, \ref{lem2}, \ref{lem3} and \ref{lem4} for 
renormalization function $R_{\infty}\left(f;p^m\right)$. For this reason we have
skipped this proof here.
\begin{theorem}\label{thm4}Let $n$ multiplicative functions $f_i(k)$, be given 
such that $f_i(k)\in{\mathbb B}\left\{G_1(N);G_2(N)\right\}$ satisfying 
(\ref{i7}). Let there exist two numbers ${\cal K}>0$ and $\gamma<a_1$ and an 
integer $r_*\geq 0$ such that $|L_r\left(\prod_{i=1}^nf_i;p^{{\bf m}}\right)|
\leq {\cal K}p^{\gamma r}$ for all $r\geq r_*$. Then
\bea
R_{\infty}\left(\prod_{i=1}^nf_i;p^{{\bf m}}\right)=\sum_{r=0}^{\infty}L_r\left(
\prod_{i=1}^nf_i;p^{{\bf m}}\right)p^{-a_1r}\;.\label{n11}
\eea
\end{theorem}
To study $R_{\infty}\left(\prod_{i=1}^nf_i;p^{{\bf m}}\right)$ in a way similar 
to the study of $R_{\infty}\left(f;p^m\right)$ in section \ref{r2} we find 
another representation for $R_{\infty}\left(\prod_{i=1}^nf_i;p^{{\bf m}}\right)$
which is different from (\ref{n11}). Substitute (\ref{n8}) into (\ref{n11}) and 
get
\bea
R_{\infty}\left(\prod_{i=1}^nf_i;p^{{\bf m}}\right)=\prod_{i=1}^nf_i\left(p^{
m_i}\right)+\left[\prod_{i=1}^nf_i\left(p^{m_i+1}\right)-L_0\left(\prod_{i=1}^n
f_i;p^{{\bf m}}\right)\prod_{i=1}^nf_i(p)\right]p^{-a_1}+\hspace{1.5cm}
\nonumber\\
\left[\prod_{i=1}^nf_i\left(p^{m_i+2}\right)-L_1\left(\prod_{i=1}^nf_i;p^{{\bf 
m}}\right)\prod_{i=1}^nf_i(p)-L_0\left(\prod_{i=1}^nf_i;p^{{\bf m}}\right)
\prod_{i=1}^nf_i(p^2)\right]p^{-2a_1}+\ldots\nonumber
\eea
Recasting the terms in the last expression we obtain
\bea
R_{\infty}\left(\prod_{i=1}^nf_i;p^{{\bf m}}\right)=\sum_{r=0}^{\infty}p^{-a_1r}
\prod_{i=1}^nf_i\left(p^{m_i+r}\right)-\sum_{r=1}^{\infty}p^{-a_1r}\prod_{i=1}^n
f_i\left(p^r\right)\cdot\sum_{r=0}^{\infty}L_r\left(\prod_{i=1}^nf_i;p^{{\bf m}}
\right)p^{-a_1r}\;.\nonumber
\eea
Thus, by comparison the last expression with formula (\ref{n11}) we get,
\bea
R_{\infty}\left(\prod_{i=1}^nf_i;p^{{\bf m}}\right)=\frac{{\bf U}\left(\prod_{i=
1}^nf_i;p^{{\bf m}}\right)}{1+{\bf V}\left(\prod_{i=1}^nf_i;p\right)}\;,
\label{n12}
\eea
where two numerical series
\bea
{\bf U}\left(\prod_{i=1}^nf_i;p^{{\bf m}}\right)=\sum_{r=0}^{\infty}
p^{-a_1r}\prod_{i=1}^nf_i\left(p^{m_i+r}\right)\;,\hspace{1cm}
{\bf V}\left(\prod_{i=1}^nf_i;p\right)=\sum_{r=1}^{\infty}p^{-a_1r}
\prod_{i=1}^nf_i\left(p^r\right)\label{n13}
\eea
are assumed to be convergent and a denominator in (\ref{n12}) does not vanish.
Formula (\ref{n12}) gives a rational representation of the renormalization
function $R_{\infty}\left(\prod_{i=1}^nf_i;p^{{\bf m}}\right)$ which is free of 
intermediate calculations of $L_r\left(\prod_{i=1}^nf_i;p^{{\bf m}}\right)$. 
\subsection{Renormalization of the Dirichlet Series $\sum_{k=1}^{\infty}f\left(
kp^m\right)k^{-s}$}\label{r42}
The Dirichlet series $D\left(f;p^m,s\right)=\sum_{k=1}^{\infty}f\left(kp^m
\right)k^{-s}$ with a scaled summation variable is a special case of summatory 
function $F\left\{f_1f_2;N,p^{{\bf m}}\right\}$, discussed in section \ref{r41},
when $f_1=f(k)$, $f_2=k^{-s}$, $m_1=m$, $m_2=0$ and $f_1f_2\in{\mathbb B}\left\{
N^0;G_2(N)\right\}$, i.e., $a_1=0$.

According to Theorem \ref{thm4} and formula (\ref{n12}) if there exist two 
numbers ${\cal K}_1>0$ and $\gamma_1<0$ and an integer $r_*\geq 0$ such that 
$|f\left(p^r\right)|\leq {\cal K}_1p^{(s+\gamma_1)r}$ for all $r\geq r_*$ then 
\bea
{\mathfrak D}\left(f;p^m,s\right)=\frac{D\left(f;p^m,s\right)}{D\left(f,s\right)
}=\sum_{r=0}^{\infty}\frac{f\left(p^{r+m}\right)}{p^{sr}}\left(1+\sum_{r=1}^{
\infty}\frac{f\left(p^{r}\right)}{p^{sr}}\right)^{-1}\;.\label{n15}
\eea
where $D\left(f,s\right)=\sum_{k=1}^{\infty}f(k)k^{-s}$ is a standard Dirichlet 
series for arithmetic function $f(k)$. Note that according to definition 
(\ref{n10}) of renormalization function $R_{\infty}\left(f\cdot k^{-s};p^m
\right)$ for the Dirichlet series the following equality holds, ${\mathfrak D}
\left(f;p^m,s\right)=p^{sm}R_{\infty}\left(f\cdot k^{-s};p^m\right)$.

We present four examples with the Dirichlet series $D\left(f;p^m,s\right)$ for 
the M\"obius $\mu(k)$, Liouville $\lambda(k)$, Euler $\phi(k)$ and divisor 
$\sigma_n(k)$ functions. Their standard Dirichlet series $D\left(f,s\right)$ 
converge to the values given in Table 2.

$\bullet\;\;$
Consider two Dirichlet series for the $\mu$-function, $D\left(\mu,s\right)=1/
\zeta(s)$ and $D\left(\mu^2,s\right)=\zeta(s)/\zeta(2s)$, $s>1$, and calculate
their scaled versions. Since $\mu\left(p^mk\right)\equiv 0$, $m\geq 2$, we 
consider here only a case $m=1$ and have $\mu^q\left(p^{r+1}\right)=(-1)^q\delta
_{0,r}$, $q=1,2$. Then
\bea
{\mathfrak D}\left(\mu;p,s\right)=-\frac{p^s}{p^s-1}\;,\quad 
{\mathfrak D}\left(\mu^2;p,s\right)=\frac{p^s}{p^s+1}\;.\label{n16}
\eea
Calculate the Dirichlet series $\Delta\left(\mu^q,s\right)=\sum_{k_1,k_2=1}
^{\infty}(k_1k_2)^{-s}\mu^q(k_1k_2)$, $q=1,2$ according to (\ref{n5}) 
\bea
\Delta\left(\mu,s\right)\!=\!\frac1{\zeta(s)}{\cal G}_2\left(\frac{\mu}{k^s}
\right),\quad {\cal G}_2\left(\frac{\mu}{k^s}\right)\!=\!\prod_{p\geq 2}\left(
1+R_{\infty}\left(\frac{\mu}{k^s};p\right)\right)\!=\!\prod_{p\geq 2}\left(1-
\frac1{p^s-1}\right)<\frac1{\zeta(s)}.\;\label{n17}
\eea
\bea
\Delta\left(\mu^2,s\right)=\frac{\zeta(s)}{\zeta(2s)}{\cal G}_2\left(\frac{
\mu^2}{k^s}\right),\quad {\cal G}_2\left(\frac{\mu^2}{k^s}\right)=\prod_{p\geq 
2}\left(1+R_{\infty}\left(\frac{\mu^2}{k^s};p\right)\right)=\prod_{p\geq 2}
\left(1+\frac1{p^s+1}\right)<\frac{\zeta(s)}{\zeta(2s)}\nonumber
\eea
that results in inequalities, $\Delta(\mu,s)<D^2(\mu,s)$ and $\Delta\left(\mu^2,
s\right)<D^2\left(\mu^2,s\right)$. A normalized product for $s=2$ in (\ref{n17})
is known as the Feller-Tornier constant $C_{FT}$ \cite{fin03},
\bea
\frac1{\zeta(2)}\prod_{p\geq 2}\left(1-\frac1{p^2-1}\right)=\prod_{p\geq 2}
\left(1-\frac{2}{p^2}\right)=C_{FT}=0.32263\;.\nonumber
\eea
$\bullet\;\;$
Consider the Dirichlet series for the $\lambda$-function, $D\left(\lambda,s
\right)=\zeta(2s)/\zeta(s)$, $s>1$, and calculate its scaled version. Keeping in
mind $\lambda\left(p^r\right)=(-1)^r$ we get ${\mathfrak D}\left(\lambda;p^m,s
\right)=(-1)^m$. Calculate the Dirichlet series $\Delta\left(\lambda,s\right)=
\sum_{k_1,k_2=1}^{\infty}(k_1k_2)^{-s}\lambda(k_1k_2)$ in accordance with 
(\ref{n5})
\bea
\Delta\left(\lambda,s\right)=\frac{\zeta(2s)}{\zeta(s)}\;{\cal G}_2\left(\frac{
\lambda}{k^s}\right),\quad {\cal G}_2\left(\frac{\lambda}{k^s}\right)=\prod_{
p\geq 2}\sum_{m=0}^{\infty}(-1)^mp^{-sm}=\prod_{p\geq 2}\frac1{1+p^{-s}}=
\frac{\zeta(2s)}{\zeta(s)}\;,\label{n17b}
\eea
i.e., $\Delta\left(\lambda,s\right)=D^2\left(\lambda,s\right)$ in accordance 
with the fact that $\lambda(k)$ is completely multiplicative.

$\bullet\;\;$
Consider the Dirichlet series for the $\phi$-function, $D\left(\phi,s\right)=
\zeta(s-1)/\zeta(s)$, $s>2$, and calculate its scaled version. Keeping in mind 
$\phi\left(p^r\right)=(p-1)p^{r-1}$ we get
\bea
{\mathfrak D}\left(\phi;p^m,s\right)=\frac{(p-1)p^{m-1}}{1-p^{-s}}\;,\quad 
m\geq 1\;.\label{n18}
\eea
$\bullet\;\;$
Consider two Dirichlet series for the $\sigma_0$-function, $D\left(\sigma_0,s
\right)=\zeta^2(s)$ and $D\left(\sigma_0^2,s\right)=\zeta^4(s)/\zeta(2s)$, $s
\geq 1$, and calculate its scaled version. Keeping in mind $\sigma_0\left(p^r
\right)=r+1$ we get
\bea
{\mathfrak D}\left(\sigma_0;p^m,s\right)&=&
(m+1)(1-p^{-s})+p^{-s}\;,\hspace{1cm}m\geq 0\;,\label{n18a}\\
{\mathfrak D}\left(\sigma_0^2;p^m,s\right)&=&
\frac{[(m+1)(1-p^{-s})+p^{-s}]^2+p^{-s}}{1+p^{-s}}\;.\nonumber
\eea
We finish this section with relationship between characteristic functions for 
multiplicative arithmetic functions $f(k)$ and $f_1(k)=f(k)\cdot k^{-s}$
\bea
L_r\left(f\cdot k^{-s};p^m\right)=L_r\left(f;p^m\right)\;p^{-(m+r)s}\;,
\label{n18b}
\eea
which follows by (\ref{b5}) and (\ref{b6}) if we substitute there the identity 
$f_1\left(p^r\right)=f\left(p^r\right)\;p^{-rs}$. 

Relation (\ref{n18b}) will be used in section \ref{r53} when calculating the 
renormalized Dirichlet series for the Ramanujan $\tau$ function.
\subsection{Renormalization of Summatory Function $\sum_{k\leq N}f\left(k^np^m
\right)$}\label{r43}
Consider the summatory function $F\left\{f,n;N,p^m\right\}=\sum_{k\leq N}f\left(
k^np^m\right)$ and derive its governing functional equation following the 
approach developed in section \ref{r12},
\bea
F\left\{f,n;N,p^m\right\}=f\left(p^m\right)\sum_{k=1\atop p\;\nmid\;k}^N
f(k^n)+\sum_{k=p\atop p\;\mid\;k}^Nf\left(p^mk^n\right)=f\left(p^m\right)\left(
\sum_{k=1}^N-\sum_{k=p\atop p\;\mid\;k}^N\right)f(k^n)+\sum_{l=1}^{N_1}f\left(
p^{m+n}l^n\right)\nonumber
\eea
which can be rewritten as follows,
\bea
F\left\{f,n;N,p^m\right\}-F\left\{f,n;N_1,p^{m+n}\right\}=f\left(p^m\right)
\left[F\left\{f,n;N,1\right\}-F\left\{f,n;N_1,p^n\right\}\right]\;.\label{n19}
\eea
In general case ($r\geq 1$) an Eq. (\ref{n19}) has a form
\bea
F\!\left\{f,n;N_r,p^{m+rn}\right\}-F\!\left\{f,n;N_{r+1},p^{m+(r+1)n}\right\}\!
=\!f\left(p^{m+rn}\right)\left[F\left\{f,n;N_r,1\right\}-F\left\{f,n;N_{r+1},
p^n\right\}\right]\nonumber
\eea
Combining last equations of running index $0\leq r\leq\left\lfloor\log_p N\right
\rfloor$ together we arrive at the functional equation,
\bea
F\left\{f,n;N_r,p^m\right\}&=&\sum_{r=0}^{\left\lfloor\log_p N\right
\rfloor}L_r\left(f,n;p^m\right)F\{f,n;N_r,1\}\;,\quad\mbox{where}\label{n20}\\
L_r\left(f,n;p^m\right)&=&f\left(p^{m+rn}\right)-\sum_{j=0}^{r-1}L_j\left(f,n;
p^m\right)f\left(p^{(r-j)n}\right).\quad\label{n21}
\eea
The straightforward calculations of $L_r\left(f;p^m\right)$ give
\bea
L_0\left(f,n;p^m\right)\!\!\!&=&\!\!\!f\left(p^m\right),\hspace{1cm}
L_1\left(f,n;p^m\right)=f\left(p^{m+n}\right)-f\left(p^m\right)f(p^n),
\label{n22}\\
L_2\left(f,n;p^m\right)\!\!\!&=&\!\!\!f\left(p^{m+2n}\right)-f\left(p^{m+n}
\right)f(p^n)-f\left(p^m\right)\left[f\left(p^{2n}\right)-f^2(p^n)\right],
\nonumber
\eea
By $n=1$ formulas (\ref{n20}), (\ref{n21})  and (\ref{n22}) are reduced to 
(\ref{b4}), (\ref{b5}) and (\ref{b6}).

Define new renormalization functions,
\bea
R\left(f,n;N,p^m\right)=\frac{F\left\{f,n;N,p^m\right\}}{F\left\{f,n;N,1\right\}
},\quad R_{\infty}\left(f,n;p^m\right):=\lim_{N\to\infty}R\left(f,n;N,p^m
\right).\quad\label{n23}
\eea
By comparison formulas (\ref{n20}), (\ref{n21}), (\ref{n22}) with (\ref{b4}),
(\ref{b5}), (\ref{b6}), respectively, and definition (\ref{n10}) with 
(\ref{n23}) we conclude that all results on renormalization of summatory 
function $F\left\{f;N,p^m\right\}$ in section \ref{r2} can be reproduced for 
summatory function $F\left\{f,n;N,p^m\right\}$ with a few necessary alterations.
Below we give (without proof) Theorem \ref{thm5} on sufficient condition to 
converge of asymptotics of renormalization function $R_{\infty}\left(f,n;p^m
\right)$. As in the case of Theorem \ref{thm4} on renormalization function $R_{
\infty}\left(f_1f_2;p^{{\bf m}}\right)$, here the proof of Theorem \ref{thm5} 
does not use new ideas and can be given following Lemmas \ref{lem1}, \ref{lem2},
\ref{lem3} and \ref{lem4} for renormalization function $R_{\infty}\left(f;p^m
\right)$. For this reason we skip it here.
\begin{theorem}\label{thm5}Let a function $f(k^n)\in{\mathbb B}\left\{G_1(N);
G_2(N)\right\}$ be given and let there exist two numbers ${\cal K}_1>0$ and 
$\gamma_1<a_1$ and an integer $r_*\geq 0$ such that $|f\left(p^r\right)|\leq 
{\cal K}_1p^{\gamma_1 r}$ for all $r\geq r_*$. Then 
\bea
R_{\infty}\left(f,n;p^m\right)=\sum_{r=0}^{\infty}\frac{L_r\left(
f,n;p^m\right)}{p^{a_1r}}\;.\label{n24}
\eea
\end{theorem}
Recasting the terms in the last expression we obtain
\bea
R_{\infty}\left(f,n;p^m\right)\!=\!\sum_{r=0}^{\infty}\frac{f\left(p^{m+nr}
\right)}{p^{a_1r}}-\sum_{r=0}^{\infty}\frac{L_r\left(f,n;p^m\right)}{p^{a_1r}}
\left[\frac{f(p^n)}{p^{a_1}}+\frac{f\left(p^{2n}\right)}{p^{2a_1}}+
\frac{f\left(p^{3n}\right)}{p^{3a_1}}+\ldots\right].\nonumber
\eea
Thus, by comparison the last expression with formula (\ref{n11}) we get,
\bea
R_{\infty}\left(f,n;p^m\right)=\frac{{\bf U}\left(f,n;p^m\right)}{1+{\bf V}
\left(f,n;p\right)}\;,\label{n25}
\eea
where two numerical series
\bea
{\bf U}\left(f,n;p^m\right)=\sum_{r=0}^{\infty}\frac{f\left(p^{m+nr}\right)}
{p^{a_1r}}\;,\hspace{1cm}
{\bf V}\left(f,n;p\right)=\sum_{r=1}^{\infty}\frac{f\left(p^{rn}\right)}
{p^{a_1r}}\;,\label{n26}
\eea
are assumed to be convergent and a denominator in (\ref{n25}) does not vanish. 

We apply formulas (\ref{n25}) and (\ref{n26}) to calculate the following 
Dirichlet series $D\left(\sigma_0,2,s,p^m\right)=\sum_{k=1}^{\infty}k^{-s}
\sigma_0(k^2p^m)$ keeping in mind \cite{apo95} the standard Dirichlet series 
$D\left(\sigma_0,2,s,1\right)=\zeta^3(s)/\zeta(2s)$, i.e., $a_1=0$. Reduce our 
problem as follows,
\bea
D\left(\sigma_0,2,s,p^m\right)=p^{sm/2}\sum_{k=1}^{\infty}\frac{\sigma_0(k^2
p^m)}{(k^2p^m)^{s/2}}=p^{sm/2}\;R_{\infty}\left(\frac{\sigma_0(k)}{k^{s/2}},2;
p^m\right)\frac{\zeta^3(s)}{\zeta(2s)}\;,\label{n27}
\eea
and calculate the renormalization function in (\ref{n27}). According to 
(\ref{n26}) we get
\bea
p^{sm/2}{\bf U}\left(\frac{\sigma_0(k)}{k^{s/2}},2;p^m\right)\!\!=\!
\frac{p^s}{p^s-1}\left(m+1+\frac{2}{p^s-1}\right),\;
1+{\bf V}\left(\frac{\sigma_0(k)}{k^{s/2}},2;p\right)\!\!=
\frac{p^s}{p^s-1}\left(1+\frac{2}{p^s-1}\right)\nonumber
\eea
so that following (\ref{n25}) and (\ref{n27}) we arrive finally at
\bea
D\left(\sigma_0,2,s,p^m\right)=\frac{(m+1)(p^s-1)+2}{p^s+1}\;\frac{\zeta^3(s)}
{\zeta(2s)}\;.\label{n28}
\eea
\section{Renormalization of the Basic Summatory Functions}\label{r5}
\vspace{-.2cm}
In this section we calculate the renormalization function $R_{\infty}\left(f;p^m
\right)$ for various summatory functions given in Tables 1, 2. For this purpose 
almost all summatory functions are treated by Theorem \ref{thm2} and 
corresponding formulas (\ref{q22}) and (\ref{q23}) based on calculation of 
$f\left(p^r\right)$. However, in section \ref{r53} we present another approach, 
which follows Theorem \ref{thm1}, and calculate the characteristic functions 
$L_r\left(f;p^m\right)$ for the Ramanujan $\tau$ function.
\subsection{Renormalization of Summatory Totient Functions}\label{r51}
\vspace{-.2cm}
In this section we apply the renormalization approach to summatory totient
functions and the Dirichlet series involving the Jordan $J_n(k)$, Euler $\phi(
k)$ and Dedekind $\psi(k)$ functions and their combinations. For the two first 
functions we make use of technical results given in \cite{fel08}. Discussing the
universality classes ${\mathbb B}\{G_1(N);G_2(N)\}$ we will skip hereafter the 
error term $G_2(N)$.
\vspace{-.5cm}
\subsubsection{Euler $\phi(k)$ function}\label{r511}
\vspace{-.2cm}
Denote the asymptotics of the summatory functions $F\left\{k^u\phi^v;N,1\right
\}=\sum_{k\leq N}k^u\phi^v(k)$ in three different ranges of varying parameters 
$-\infty<u<\infty$ and $v\in{\mathbb Z}$ which is given in \cite{fel08},
\bea
F\left\{k^u\phi^v;N,1\right\}\stackrel{N\to\infty}{=}\left\{\begin{array}
{lrr}A(u,v)N^{u+v+1}&,&u+v>-1\;,\\B(u,v)\ln N&,&u+v=-1\;,\\
C(u,v)&,&u+v<-1\;,\end{array}\right.\label{r512a}
\eea
such that $A(u,0)=(u+1)^{-1}$, $B(-1,0)=1$ and $C(u,0)=\zeta(-u)$. 

In the case 
of arbitrary $u$ and $v$ a lot of expressions for $A(u,v)$, $B(u,v)$ and 
$C(u,v)$ can be found in \cite{apo95}, \cite{chl27}, \cite{now89}, \cite{san06} 
and \cite{sur82}. Here we focus on renormalization functions not specifying the 
explicit expressions. According to \cite{fel08}, if $v\neq 0$ then $A(u,v)$, 
$B(u,v)$ and $C(u,v)$ are bounded from above as follows: if $v\in{\mathbb Z}_+$,
then
\bea
0<A(u,v)\leq\frac{(u+v+1)^{-1}}{\zeta(v+1)},\quad 0<B(u,v)\leq 1/\zeta(v+1),
\quad 0<C(u,v)\leq\zeta(-u-v)/\zeta(-u)\;,\nonumber
\eea
and if $v\in{\mathbb Z}_-$, then
\bea
&&(u+v+1)^{-1}<A(u,v)\leq 2^{\frac{|v|}{2}}\;{\cal D}_{\infty}(v,1)(u+v+1)^{-1},
\quad 1<B(u,v)\leq 2^{\frac{|v|}{2}}\;{\cal D}_{\infty}(v,1)\;,\nonumber\\
&&\zeta(-u-v)<C(u,v)\leq2^{\frac{|v|}{2}}\;{\cal D}_{\infty}(v,-u-v)\zeta(-u-v)
\;,\quad\mbox{where}\quad {\cal D}_{\infty}(v,s)=\prod_{r=1}^{|v|}\zeta\left(s+
\frac{r}{2}\right).\quad\nonumber
\eea
Calculate their renormalization functions $R_{\infty}\left(k^u\phi^v;p^m\right)$
according to (\ref{q22}) and (\ref{q23}),
\bea
R_{\infty}\left(k^u\phi^v;p^m\right)&=&p^{m(u+v)+1}\frac{(p-1)^{v-1}}{p^v+(p-1)
^{v-1}}\;,\quad u+v\geq -1\;,\label{r512b}\\
R_{\infty}\left(k^u\phi^v;p^m\right)&=&p^{(m-1)(u+v)}\frac{(p-1)^v}{p^{-u}-p^v+
(p-1)^v}\;,\quad u+v<-1\;.\nonumber
\eea
By (\ref{r512b}) and (\ref{r511b}) we have an equality $R_{\infty}\left(k^uJ_1;
p^m\right)=R_{\infty}\left(k^u\phi;p^m\right)$ in all ranges of $u$.
\subsubsection{Jordan $J_v(k)$ and Dedekind $\psi(k)$ functions}\label{r512}
Regarding the Jordan function $J_v(k)$, asymptotics of the summatory functions 
$F\left\{k^uJ_v;N,1\right\}=\sum_{k\leq N}k^uJ_v(k)$, $v\in{\mathbb Z}_+$, can 
be given in three different ranges of varying parameters \cite{fel08},
\bea
F\left\{k^uJ_v;N,1\right\}\stackrel{N\to\infty}{\simeq}\left\{\begin{array}{lrr}
(u+v+1)^{-1}\frac{N^{u+v+1}}{\zeta(v+1)}&,&u+v>-1\;,\\
\ln N/\zeta(v+1)&,&u+v=-1\;,\\
\zeta(-u-v)/\zeta(-u)&,&u+v<-1\;,\end{array}\right.\label{r511a}
\eea
and calculate their renormalization function $R_{\infty}\left(k^uJ_v;p^m\right)$
in accordance with (\ref{q22}) and (\ref{q23})
\bea
\frac{R_{\infty}\left(k^uJ_v;p^m\right)}{p^v-1}=\frac{p^{m(u+v)+1}}{p^{v+1}-1},
\;u+v\geq 
-1\;,\quad\frac{R_{\infty}\left(k^uJ_v;p^m\right)}{p^v-1}=\frac{p^{(
m-1)(u+v)}}{p^{-u}-1},\;u+v<-1\;.\quad\label{r511b}
\eea
Consider another summatory functions $F\left\{k^u\psi^v;N,1\right\}=\sum_{k\leq 
N}k^u\psi^v(k)$ in three different ranges of varying parameters $-\infty<u<
\infty$ and $v\in{\mathbb Z}$ and note that
\bea
k^u\psi^v(k)\in\left\{\begin{array}{lrr}
{\mathbb B}\{N^{u+v+1}\}&,&u+v>-1\;,\\{\mathbb B}\{\ln N\}&,&u+v=-1\;,\\
{\mathbb B}\{N^0\}&,&u+v<-1\;.\end{array}\right.\label{r513a}
\eea
The explicit asymptotics for some summatory functions $F\left\{k^u\psi^v;N,1
\right\}$ are given in \cite{sita79}, \cite{sur82}. 

Calculate their renormalization functions $R_{\infty}\left(k^u\psi^v;p^m\right)$ according to 
(\ref{q22}) and (\ref{q23}),
\bea
R_{\infty}\left(k^u\psi^v;p^m\right)&=&p^{m(u+v)+1}\frac{(p+1)^v}{(p+1)^v+(p-1)
p^v}\;,\quad u+v\geq -1\;,\label{r513b}\\
R_{\infty}\left(k^u\psi^v;p^m\right)&=&p^{(m-1)(u+v)}\frac{(p+1)^v}{(p+1)^v-
p^v+p^{-u}}\;,\quad u+v<-1\;.\nonumber
\eea
We finish this section with summatory $F\left\{(\phi/\psi)^v;N,1\right\}=
\sum_{k\leq N}\phi^v(k)\psi^{-v}(k)$, $v\in{\mathbb Z}$. 

The asymptotics $N^{-1}F\left\{\phi/\psi;N,1\right\}\simeq\prod_{p}(1-2/(p(p+1))
\simeq 0.4716$ is known due to \cite{sur82}. Keeping in mind $(\phi/\psi)^v\in
{\mathbb B}\{N\}$ calculate a corresponding renormalization function,
\bea
R_{\infty}\left(\left(\frac{\phi}{\psi}\right)^v;p^m\right)=
\frac{p(p-1)^{v-1}}{(p-1)^{v-1}+(p+1)^v}\;,\label{r513c}
\eea
which does not dependent on $m$.
\subsection{Renormalization of Summatory Non-Totient Functions}\label{r52}
In this section we apply the renormalization approach to summatory functions
and Dirichlet series involving divisor $\sigma_a(k)$, prime divisor $\beta(k)$,
Piltz $d_n(k)$, abelian group enumeration $\alpha(k)$ functions, Ramanujan sum 
$C_q(n)$ and some of their their combinations.
\subsubsection{Divisor function $\sigma_a(k)$ and prime divisor function 
$\beta(k)$}\label{r521}
The divisor function $\sigma_a(k)$ is defined as a sum of the $a$th powers of 
the divisors of $k$. For $k=p^r$ we have $\sigma_a\left(p^r\right)=\left(p^{a(
r+1)}-1\right)/\left(p^a-1\right)$, $a\neq 0$, and $\sigma_0\left(p^r\right)=
r+1$.

$\bullet\;\;$ $F\left\{k^{-s}\sigma_a;N,1\right\}$, $a>0$, $s\geq 1+a$, 
$k^{-s}\sigma_a\in{\mathbb B}\left\{N^0\right\}$
\bea
R_{\infty}\left(\frac{\sigma_a}{k^s};p^m\right)=\frac{p^{a(m+1)+s}-p^{a(m+1)}-
p^s+p^a}{p^{(m+1)s}(p^a-1)}\;.\label{r521a}
\eea
$\bullet\;\;$ $F\left\{\sigma_a;N,1\right\}$, $a>0$, $\sigma_a\in{\mathbb B}
\left\{N^{a+1}\right\}$
\bea
R_{\infty}\left(\sigma_a;p^m\right)=\frac{p^{a(m+1)+1}-p^{am}-p+1}{p\;(p^a-1)}
\;.\label{r521b}
\eea
$\bullet\;\;$ $F\left\{\sigma_a;N,1\right\}$, $a<0$, $\sigma_a\in{\mathbb B}
\{N\}$
\bea
R_{\infty}\left(\sigma_a;p^m\right)=\frac{p^{am+1}-p^{am}-p^{1-a}+1}{p^{1-a}\;
(p^a-1)}\;.\label{r521c}
\eea
$\bullet\;\;$ $F\left\{\sigma_0^n;N,1\right\}$, $\sigma_0^n\in{\mathbb B}\left\{
N\left(\log N\right)^{2^n-1}\right\}$
\bea
R_{\infty}\left(\sigma_0^n;p^m\right)=\frac{{\cal S}(n,p,m+1)}{{\cal S}(n,p,1)}
\;,\quad {\cal S}(n,p,t)=\sum_{k=1}^n{n\choose k}t^{n-k}Li_{-k}(p^{-1})+
\frac{t^n}{1-p^{-1}}\;,\label{r521d}
\eea
where $Li_{-k}(x)$ is the polylogarithm function defined in (\ref{b11a1}). 
Substituting $n=1$ into (\ref{r521d}) we get $R_{\infty}\left(\sigma_0;p^m
\right)=m+1-m/p$. According to (\ref{r521b}) and (\ref{r521c}), this expression 
coincides with both limits of $R_{\infty}\left(\sigma_a;p^m\right)$, when $a\to 
0$, for $\sigma_a>$ and $\sigma_a<0$, respectively.

$\bullet\;\;$ $F\left\{\sigma_0\sigma_a;N,1\right\}$, $a>0$, $\sigma_0\sigma_a
\in{\mathbb B}\left\{N^{a+1}\log N\right\}$
\bea
R_{\infty}\left(\sigma_0\sigma_a;p^m\right)\!=\!\frac{p^{am+1}\left(p^{a+1}-1
\right)^2\!-p^{a+1}(p-1)^2\!+m(p-1)\left(p^{a+1}-1\right)\left(p^{am}\left(p^{a+
1}-1\right)\!-p+1\right)}{p\left(p^a-1\right)\left(p^{a+2}-1\right)}\nonumber
\eea
such that $R_{\infty}\left(\sigma_0\sigma_a;1\right)=1$. Note that $R_{\infty}
\left(\sigma_0\sigma_a;p^m\right)\stackrel{a\to 0}{\to}R_{\infty}\left(\sigma_0
^2;p^m\right)$ according to (\ref{r521d}).

$\bullet\;\;$ $F\left\{\sigma_a^2;N,1\right\}$, $a>0$, $\sigma_a^2\in{\mathbb B}
\left\{N^{2a+1}\right\}$
\bea
R_{\infty}\left(\sigma_a^2;p^m\right)=\frac{(p-1)\left(p^{a+1}-1\right)+p^{am}
\left(p^{1+2 a}-1\right)\left(p^{a m}\left(p^{1+a}-1\right)-2(p-1)\right)}
{p\left(p^a-1\right)^2\left(p^{a+1}+1\right)}\;,\label{r521j}  
\eea
that gives $R_{\infty}\left(\sigma_a^2;p^m\right)\stackrel{a\to 0}{\to}R_{
\infty}\left(\sigma_0^2;p^m\right)=(p+(m(p-1)+p)^2)/(p(p+1))$ according to 
(\ref{r521d}).

$\bullet\;\;$ $F\left\{1/\sigma_0;N,1\right\}$, $1/\sigma_0\in{\mathbb B}\left\{
N/\sqrt{\log N}\right\}$
\bea
R_{\infty}\left(\frac1{\sigma_0};p^m\right)=-\frac{\;_2F_1\left(m+1,1;m+2;p^{-1}
\right)}{(m+1)\;p\;\ln\left(1-p^{-1}\right)}\;,\label{r521e}
\eea
where $_2F_1(a,b;c;z)$ denotes a generalized hypergeometric function 
\cite{em81}. By (\ref{r521e}) we get for $m=1$
\bea
R_{\infty}\left(\frac1{\sigma_0};p\right)=p+\frac1{\ln\left(1-p^{-1}\right)}
\;,\quad R_{\infty}\left(\frac1{\sigma_0};p\right)\stackrel{p\to\infty}
{\longrightarrow}\;\frac1{2}\;,\quad R_{\infty}\left(\frac1{\sigma_0};p\right)
\stackrel{p\to 1}{\longrightarrow}\;1\;.\label{r521f}
\eea
Expanding an expression (\ref{r521e}) as an infinite series $\sum_{r=0}^{\infty}
L_r\left(1/\sigma_0;p^m\right)p^{-r}$ in accordance with Theorem \ref{thm1}, one
can calculate $L_r\left(1/\sigma_0;p^m\right)$ and verify that for $0\leq r\leq 
3$ they coincide with those given in (\ref{b6b}).
 
$\bullet\;\;$ $F\left\{\sigma_1/\sigma_0;N,1\right\}$, $\sigma_1/\sigma_0\in
{\mathbb B}\left\{N^2/\sqrt{\log N}\right\}$
\bea
R_{\infty}\left(\frac{\sigma_1}{\sigma_0};p^m\right)=\frac{p^{m+1}\;_2F_1\left(
m+1,1;m+2;p^{-1}\right)-\;_2F_1\left(m+1,1;m+2;p^{-2}\right)}{(m+1)\;p^2\;
\log\left(1+p^{-1}\right)}\;,\label{r521g}
\eea
that gives for $m=1$
\bea
R_{\infty}\left(\frac{\sigma_1}{\sigma_0};p\right)=p^2-\frac{p-1}{\log\left(1+
p^{-1}\right)}\;,\quad \frac1{p}R_{\infty}\left(\frac{\sigma_1}{\sigma_0};
p\right)\stackrel{p\to\infty}{\longrightarrow}\;\frac1{2}\;.\label{r521h}
\eea
$\bullet\;\;$ $F\left\{\sigma_1/\phi;N,1\right\}$, $F\left\{\sigma_1/\psi;N,1
\right\}$, $\quad\sigma_1/\phi,\;\sigma_1/\psi\in{\mathbb B}\left\{N\right\}$,
$m\geq 1$.
\bea
R_{\infty}\left(\frac{\sigma_1}{\phi};p^m\right)=\frac{p^3\left(1+p-p^{-m}
\right)}{p^4-p^3 +p^2+ p-1}\;,\quad 
R_{\infty}\left(\frac{\sigma_1}{\psi};p^m\right)=\frac{p^3\left(1+p-p^{-m}
\right)}{p^4+p^3-p^2-p+1}\;.\label{r521j}
\eea

The prime divisor function $\beta(k)$ is defined by formula $\beta\left(p^{a_1}
\cdot\ldots\cdot p^{a_n}\right)=a_1\cdot\ldots\cdot a_n$.

$\bullet\;\;$ $F\left\{\beta;N,1\right\}$, $\beta\in{\mathbb B}\left\{N\right\}$
\bea
R_{\infty}\left(\beta;p^m\right)=p\;\frac{m(p-1)+1}{p^2-p+1}\;.\label{r524a}
\eea
Let us note a curious consequence of (\ref{r524a}) when $m=p$ : $R_{\infty}
\left(\beta;p^p\right)=p$.

$\bullet\;\;$ $F\left\{k^{-s}\beta;N,1\right\}$, $\beta\in{\mathbb B}\left\{   
N^0\right\}$
\bea
R_{\infty}\left(\frac{\beta}{k^s};p^m\right)=p^{(1-m)s}\;\frac{m(p^s-1)+1}
{p^{2s}-p^s+1}\;.\label{r524b}
\eea

A generalized summatory function $F\left\{\beta^n;N,1\right\}$ was considered in
\cite{kno73} such that ${\cal F}(\beta^n)=\prod_{p\geq 2}\left(1+\sum_{j=2}^{
\infty}\left[j^n-(j-1)^k\right]p^{-j}\right)$, e.g., ${\cal F}(\beta)=\prod_{
p\geq 2}\left[1-1/(p(p-1))\right]=\zeta(2)\zeta(3)/\zeta(6)$.   

$\bullet\;\;$ $F\left\{\beta^n;N,1\right\}$, $\beta^n\in{\mathbb B}\left\{N
\right\}$
\bea
R_{\infty}\left(\beta^n;p^m\right)=\frac{{\cal S}(n,p,m)}{1+Li_{-n}(p^{-1})}\;,
\label{r524c}
\eea
where ${\cal S}(n,p,m)$ is defined in (\ref{r521d}).
\subsubsection{Piltz function $d_n(k)$ and the sum of two squares function 
$r_2(k)$}\label{r522}
The Piltz function $d_n(k)$ is defined as a number of ways to write the positive
integer $k$ as a product of $n$ (positive integer) factors. For $k=p^r$ we have 
$d_n\left(p^r\right)={n+r-1\choose r}$. By definition, it holds $d_1(k)=1$ and 
$d_2(k)=\sigma_0(k)$.

$\bullet\;\;$ $F\left\{k^{-s}d_n;N,1\right\}$, $d_n\in{\mathbb B}\left\{N^0
\right\}$
\bea
R_{\infty}\left(\frac{d_n}{k^s};p^m\right)=p^{-sm}\;{n+m-1\choose m}\;_2F_1
\left(m,1-n;m+1;p^{-s}\right)\;.\label{r522a}
\eea
$\bullet\;\;$ $F\left\{d_n;N,1\right\}$, $d_n\in{\mathbb B}\left\{N\left(\log 
N\right)^{n-1}\right\}$
\bea
R_{\infty}\left(d_n;p^m\right)={n+m-1\choose m}\;_2F_1\left(m,1-n;m+1;p^{-1}
\right)\;.\label{r522b}
\eea
Note that $R_{\infty}\left(d_2;p^m\right)$ is coincided with $R_{\infty}\left(
\sigma_0;p^m\right)$ given in section \ref{r521}.

$\bullet\;\;$ $F\left\{d_n^2;N,1\right\}$, $d_n^2\in{\mathbb B}\left\{N\left(
\log N\right)^{n^2-1}\right\}$
\bea
R_{\infty}\left(d_n^2;p^m\right)={n+m-1\choose m}^2\frac{\;_3F_2\left(1,m+n,
m+n;m+1,m+1;p^{-1}\right)}{\;_2F_1\left(n,n;1;p^{-1}\right)}\;,\label{r522c}
\eea
and $R_{\infty}\left(d_2^2;p^m\right)=R_{\infty}\left(\sigma_0^2;p^m\right)$ in 
accordance with (\ref{r521d}).

$\bullet\;\;$ $F\left\{1/d_n;N,1\right\}$, $1/d_n\in{\mathbb B}\left\{N\left(
\log N\right)^{1/n-1}\right\}$
\bea
R_{\infty}\left(\frac1{d_n};p^m\right)=np\;{n+m-1\choose m}^{-1}\frac{\;_2F_1
\left(1,m+1;m+n;p^{-1}\right)}{n\;p+\;_2F_1\left(2,1;n+1;p^{-1}\right)}\;,
\label{r522d}
\eea
such that $R_{\infty}\left(d_2^{-1};p^m\right)=R_{\infty}\left(\sigma_0^{-1};
p^m\right)$ in accordance with (\ref{r521e}).

The number of representations of $k$ by two squares, allowing zeros and 
distinguishing signs and order, is denoted by $r_2(k)$. If $k=p^r$ then,
\bea
r_2\left(p^r\right)=\left\{\begin{array}{lrl}4(r+1)&,&p=1\pmod 4\;,\\
4&,&p=2\end{array}\right.\;\quad\mbox{and}\quad r_2\left(p^r\right)
=\left\{\begin{array}{lrl}4&,&p=3\pmod 4\;,\;2\mid r\\0&,&p=3\pmod 4\;,\;2\nmid 
r\end{array}\right.\nonumber
\eea

$\bullet\;\;$ $F\left\{r_2;N,1\right\}$, $r_2\in{\mathbb B}\left\{N\right\}$
\bea
R_{\infty}\left(r_2;p^m\right)\!=\!\left\{\begin{array}{lrl}\!m+1-m/p\!\!\!&,&\!
\!p=1\!\!\!\pmod 4\;,\\\!1\!\!\!&,&\!\!p=2\;,\end{array}\right.
R_{\infty}\left(r_2;p^m\right)\!=\!\left\{\begin{array}{lrl}\!1\!\!&,&\!\!p=3
\!\!\!\pmod 4,\;2\mid m\\\!1/p\!\!&,&\!\!p=3\!\!\!\pmod 4,\;2\nmid m
\end{array}\right.\nonumber
\eea
$\bullet\;\;$ $F\left\{r_2/k;N,1\right\}$, $r_2/k\in{\mathbb B}\left\{\log N
\right\}$
\bea
R_{\infty}\left(\frac{r_2}{k};p^m\right)=p^{-m}\;R_{\infty}\left(r_2;p^m\right)
\label{r523a}
\eea
\subsubsection{Ramanujan sum $C_q(n)$ and Abelian group enumeration function 
$\alpha(k)$}\label{r523}
Ramanujan's sum $C_q(n)$, $q,n\geq 1$, is a multiplicative arithmetic function 
which is defined by formula $C_q(n)=\sum_{a=1}^q\exp(2\pi ina/q)$, $(a,q)=1$ 
such that $C_{q_1q_2}(n)=C_{q_1}(n)C_{q_2}(n)$ if $(q_1,q_2)=1$ and $C_{1}(n)=
1$. If $q=p^r$ then
\bea
C_{p^r}(n)=0,\;\;\mbox{if}\;\;p^{r-1}\nmid n;\quad C_{p^r}(n)=-p^{r-1},\;\;
\mbox{if}\;\;p^{r-1}\mid n,\;p^r\nmid n;\quad C_{p^r}(n)=\phi(r),\;\;\mbox{if}
\;\;p^r\mid n,\nonumber\\
C_{p}(n)=-1,\;\;\mbox{if}\;\;p\nmid n;\quad C_{p^r}(n)=0,\;\;\mbox{if}\;\;
p\nmid n\;\;\mbox{and}\;\;r\geq 2\;.\hspace{5.5cm}\label{r525a}
\eea
Consider the Dirichlet summatory function for $C_q(n)$ : $F\left\{C_q(n)q^{-s};
\infty,1,s\right\}=\sum_{q=1}^{\infty}C_q(n)q^{-s}$, $s>1$, where $n$ is kept 
constant. It is convergent to $n^{-s+1}\sigma_{s-1}(n)\zeta^{-1}(s)$. We find a 
rescaled summatory function $F\left\{C_{q}(n)q^{-s};\infty,p^m,s\right\}=
\sum_{q=1}^{\infty}C_{p^mq}(n)q^{-s}$.
 
Let a number $n$ is such that $p^a\mid n$ but $p^{a+1}\nmid n$, $a\in{\mathbb Z}
_+$, and $a\geq m\geq 1$, then
\bea
{\bf U}\left(\frac{C_q(n)}{q^s};p^m\right)\!=\!p^{-ms}\left(\sum_{r=m}^a\frac{
\phi(p^r)}{p^{sr}}-\frac{p^a}{p^{s(a+1)}}\right),\;\;{\bf V}\left(\frac{C_q(n)}
{q^s};p^m\right)\!=\!\sum_{r=1}^a\frac{\phi(p^r)}{p^{sr}}-\frac{p^a}{p^{s(a+1)}}
\quad\label{r525y}
\eea
that gives
\bea
F\left\{\frac{C_q(n)}{q^s};\infty,p^m,s\right\}=\frac{p^s-1+p^{(s-1)(a+2-m)}-1)}
{(p^s-1)(p^{(s-1)(a+1})}\;F\left\{\frac{C_q(n)}{q^s};\infty,1,s\right\}\;.
\label{r525x}
\eea
If $a\geq m=0$ then by (\ref{r525y}) we have ${\bf U}\left(C_q(n)q^{-s};1
\right)=1+{\bf V}\left(C_q(n)q^{-s};1\right)$, and therefore $F\left\{C_{q}(n)
q^{-s};\infty,p^m,s\right\}\stackrel{m\to 0}{\to}F\left\{C_q(n)q^{-s};\infty,
1,s\right\}$.

If $m\geq 2$ and a number $n$ is not divided by $p$, i.e., $a=0$, then $F\left\{
C_{q}(n)q^{-s};\infty,p^m,s\right\}=0$. 

Finally, if $m=1$ and $a=0$, then $F\left\{C_q(n)q^{-s};\infty,p,s\right\}=
F\left\{C_q(n)q^{-s};\infty,1,s\right\}/\left(p^{-s}-1\right)$.

Abelian group enumeration function $\alpha(k)$ accounts for the number of
(isomorphism classes of) commutative groups of order $k$. By definition it 
satisfies $\alpha\left(p^r\right)={\cal P}(r)$, where ${\cal P}(r)$ denotes an
unrestricted partition function \cite{and76} and ${\cal P}(0)=1$.

$\bullet\;\;$ $F\left\{\alpha;N,1\right\}$, $\alpha\in{\mathbb B}\left\{N
\right\}$
\bea
R_{\infty}\left(\alpha;p^m\right)=p^m\left(1-Q\left(p^{-1}\right)\sum_{k=0}^{
m-1}{\cal P}(k)p^{-k}\right),\quad Q(x)=\prod_{j=1}^{\infty}\left(1-x^j\right)
\quad |x|\leq 1\;.\label{r524d}
\eea
$\bullet\;\;$ $F\left\{1/\alpha;N,1\right\}$, $1/\alpha\in{\mathbb B}\left\{N
\right\}$
\bea
R_{\infty}\left(\frac1{\alpha};p^m\right)=p^m\;\frac{{\cal T}(p,m)}{{\cal T}
(p,0)},\quad {\cal T}(p,t)=\sum_{k=t}^{\infty}\frac{p^{-k}}{{\cal P}(k)}\;.
\label{r524e}
\eea
\subsection{Renormalization of Summatories Associated with Ramanujan's $\tau$ 
Function}\label{r53}
The Ramanujan $\tau$ function is a multiplicative arithmetic function which is 
mostly known due to its generating function, $\sum_{k=1}^{\infty}\tau(k)x^k=x
\prod_{k=1}^{\infty}(1-x^k)^{24}$, $|x|<1$. For our purpose to calculate the 
renormalization function for any summatory function $F\{f(\tau,k);N,1\}$ with 
$f(\tau,k)$ involving the $\tau$ function, it is important to know its 
recursive relation for $k=p^r$,
\bea
\tau\left(p^r\right)=\sum_{j=0}^{\left\lfloor r/2\right\rfloor}(-1)^j{r-j\choose
r-2j}p^{11j}\tau^{r-2j}(p)\;.\label{r531}
\eea
Then, making use of formulas (\ref{q23}) and (\ref{q24}) we can arrive at $R_{
\infty}\left(f(\tau,k);p^m\right)$ due to the finite computational procedure.
However, the representation (\ref{r531}) is too difficult to make worth, so we
choose another way to find $R_{\infty}\left(f(\tau,k);p^m\right)$, namely, by 
Theorem \ref{thm1} and by calculating the characteristic functions $L_r\left(
\tau;p^m\right)$. Start with identity for $\tau$ function \cite{leh43}
\bea
\tau\left(p^{r+1}\right)=\tau(p)\tau\left(p^r\right)-p^{11}\tau\left(p^{r-1}
\right)\;.\label{r532}
\eea
The following Proposition is based on recursion (\ref{b5}) and the last 
identity.
\begin{proposition}\label{pro2}
\bea
L_0\left(\tau;p^m\right)=\tau\left(p^m\right),\quad L_1\left(\tau;p^m\right)=-
p^{11}\tau\left(p^{m-1}\right),\quad L_r\left(\tau;p^m\right)=0,\quad r\geq 2\;.
\label{r533}
\eea
\end{proposition}
{\sf Proof} $\;\;\;$Calculating the four first expressions of $L_r\left(\tau;p^m
\right)$ one can verify that (\ref{r533}) holds, i.e., $L_2\left(\tau;p^m\right)
=L_3\left(\tau;p^m\right)=0$. Prove by induction that $L_r\left(\tau;p^m\right)
=0$, $r\geq 2$. 

Indeed, let $L_q\left(\tau;p^m\right)=0$ for $2\leq q\leq r$, then keeping in 
mind (\ref{b5}) and (\ref{r532}) write this equality in another representation,
\bea
L_q\left(\tau;p^m\right)=\tau\left(p^{m+q}\right)-\tau\left(p^m\right)\tau\left(
p^q\right)+p^{11}\tau\left(p^{m-1}\right)\tau\left(p^{q-1}\right)=0,\quad 2\leq 
q\leq r\;.\label{r534}
\eea
Making use of (\ref{r532}) and (\ref{r534}) write the next term $L_{r+1}\left(
\tau;p^m\right)$,
\bea
L_{r+1}\left(\tau;p^m\right)=\tau\left(p^{m+r+1}\right)-\tau\left(p^m\right)
\tau\left(p^{r+1}\right)+p^{11}\tau\left(p^{m-1}\right)\tau\left(p^r\right)\;,
\nonumber
\eea
and calculate a difference,
\bea
L_{r+1}\left(\tau;p^m\right)\!-\!\tau(p)L_r\left(\tau;p^m\right)\!=\!\left[
\tau\left(p^{m+r+1}\right)-\tau(p)\tau\left(p^{m+r}\right)\right]\!-\!\tau\left(
p^m\right)\left[\tau\left(p^{r+1}\right)-\tau(p)\tau\left(p^{r}\right)\right]
\nonumber\\
+p^{11}\tau\left(p^{m-1}\right)\left[\tau\left(p^r\right)-
\tau(p)\tau\left(p^{r-1}\right)\right]\;.\hspace{3cm}\label{r536}
\eea
By identity (\ref{r532}) the r.h.s. in equality (\ref{r536}) can be 
reduced as follows,
\bea
L_{r+1}\left(\tau;p^m\right)\!-\!\tau(p)L_r\left(\tau;p^m\right)\!=-p^{11}\left[
\tau\left(p^{m+r-1}\right)-\tau\left(p^m\right)\tau\left(p^{r-1}\right)-p^{11}
\tau\left(p^{m-1}\right)\tau\left(p^{r-2}\right)\right]\;.\nonumber
\eea
By comparison with (\ref{r534}) one can recognize the function $L_{r-1}\left(
\tau;p^m\right)$ in the brackets of the last expression. Then, combining this 
fact with (\ref{r536}) and assumption (\ref{r534}) we get
\bea
L_{r+1}\left(\tau;p^m\right)=\tau(p)L_r\left(\tau;p^m\right)-p^{11}L_{r-1}\left(
\tau;p^m\right)=0\;.\label{r537}
\eea
Thus, proof is finished.$\;\;\;\;\;\;\Box$

An important spinoff arising from Proposition \ref{pro2} is the identity 
(\ref{r534}) which does generalize an identity (\ref{r532}).
\begin{corollary}\label{cor2}
The Ramanujan $\tau$ function satisfies an identity,
\bea
\tau\left(p^{m+n}\right)=\tau\left(p^m\right)\tau\left(p^n\right)-p^{11}\tau
\left(p^{m-1}\right)\tau\left(p^{n-1}\right),\quad m,n\geq 1\;.\label{r538}
\eea
\end{corollary}
The last statement implies two inequalities which could easily be veryfied. The 
1st inequality looks quite trivial, $\tau\left(p^{2n}\right)<\tau^2\left(p^n
\right)$. Regarding the 2nd inequality, let $p_{\ast}$ and $n_{\ast}$ be choosen
in such a way that $\tau\left(p_{\ast}^{2n_{\ast}}\right)<0$, e.g., $\tau\left(
2^2\right),\tau\left(3^2\right),\tau\left(5^2\right),\tau\left(7^2\right)<0$, 
$\;\tau\left(5^4\right),\tau\left(11^4\right)<0$ etc. Then the following 
inequality holds,
\bea
|\tau\left(p_{\ast}^{2n_{\ast}}\right)|<p_{\ast}^{11}\tau^2\left(p_{\ast}^{
n_{\ast}-1}\right)\;.\label{r538x}
\eea
In the case $n_{\ast}=2$ let us combine (\ref{r538x}) with Deligne's inequality 
$|\tau(p)|<2p^{11/2}$ for the Ramanujan $\tau$ functions \cite{del74} and get, 
\bea
|\tau\left(p_{\ast}^4\right)|<4p_{\ast}^{22}\;.\label{rd1}
\eea
Note that the upper bound in (\ref{rd1}) is stronger than the bound which came 
by combining (\ref{r532}) and Deligne's inequality for arbitrary prime $p$. 
Indeed, 
\bea
\tau\left(p^4\right)=\tau^4(p)-3p^{11}\tau^2(p)+p^{22}<17p^{22}-3p^{11}\tau^2(p)
\;.\label{rd2}
\eea
However, $4p^{22}<17p^{22}-3p^{11}\tau^2(p)$ that follows by Deligne's 
inequality, $17p^{22}-3p^{11}\tau^2(p)>5p^{22}$.

Straightforward numerical calculations show that the inequality (\ref{rd1}) 
holds also for the first 474 primes irrespectively whether the requirement 
$\tau\left(p^4\right)<0$ is holding,
\bea
p_{475}=3371\;,\quad \frac{\tau\left(3371^4\right)}{4\cdot 3371^{22}}\simeq 
1.0119\;.\nonumber
\eea
\subsubsection{Renormalization of the Ramanujan $\tau$ Dirichlet series}
\label{rr531}
Recalling the Ramanujan conjecture on $\tau$ function, $\tau(N)={\cal O}\left(
N^{11/2+\varepsilon}\right)$ proved by Deligne, consider the $\tau$ Dirichlet 
series $F\left\{\tau\cdot k^{-s};\infty,p^m\right\}=\sum_{k=1}^{\infty}\tau(k)
\;k^{-s}$, $\tau\cdot k^{-s}\in{\mathbb B}\left\{N^0\right\}$, which converges 
absolutely if $s>13/2$. By Theorem \ref{thm1} and relationship (\ref{n18b}) 
between characteristic functions for $\tau(k)$ and $\tau(k)\;k^{-s}$ the 
renormalization function $R_{\infty}\left(\tau\cdot k^{-s};p^m\right)$ reads
\bea
R_{\infty}\left(\tau\cdot k^{-s},p^m\right)=\sum_{r=0}^{\infty}L_r\left(\tau
\cdot k^{-s};p^m\right)=\sum_{r=0}^{\infty}L_r\left(\tau;p^m\right)p^{-s(r+m)}
\;.\label{r538a}
\eea
Applying Proposition \ref{pro2} to (\ref{r538a}) we obtain,
\bea
R_{\infty}\left(\tau\cdot k^{-s},p^m\right)=p^{-sm}\left(\tau\left(p^m\right)-
\frac{\tau\left(p^{m-1}\right)}{p^{s-11}}\right)\;.\nonumber
\eea
Making use of relation between renormalization function $R_{\infty}\left(f\cdot 
k^{-s};p^m\right)$ and a ratio ${\mathfrak D}\left(f;p^m,s\right)$ between 
scaled and unscaled Dirichlet series, given in section \ref{r42}, we get 
finally,
\bea
\sum_{k=1}^{\infty}\frac{\tau(kp^m)}{k^s}=\left(\tau\left(p^m\right)-\frac{
\tau\left(p^{m-1}\right)}{p^{s-11}}\right)\sum_{k=1}^{\infty}\frac{\tau(k)}{
k^s}\;.\label{r539}
\eea
Formula (\ref{r539}) gives rise to several special cases, e.g.,
\bea
\sum_{k=1}^{\infty}\frac{\tau(kp)}{k^s}=\left(\tau(p)-p^{11-s}\right)
\sum_{k=1}^{\infty}\frac{\tau(k)}{k^s}\;,\hspace{1cm}\sum_{n=0}^m\sum_{k=1}^{
\infty}\frac{\tau(kp^n)}{k^{11}}=\tau\left(p^m\right)\sum_{k=1}^{\infty}
\frac{\tau(k)}{k^{11}}\;.\nonumber
\eea
\subsubsection{Renormalization functions $R_{\infty}\left(\tau^2,p^m\right)$
and $R_{\infty}\left(\tau^2\cdot k^{-25/2},p^m\right)$}\label{rr532}
In this section we renormalize the summatory functions $F\left\{\tau^2;N,p^m
\right\}$ and $F\left\{\tau^2\cdot k^{-25/2};N,p^m\right\}$ given in Tables 1 
and 2. For this purpose we start by calculating the characteristic functions 
$L_r\left(\tau^2;p^m\right)$. In contrast to $L_r\left(\tau;p^m\right)$ 
described in Proposition \ref{pro2} the present case is not so simple but still
allows to find the general expressions.
\begin{proposition}\label{pro3}
\bea
&&L_0\left(\tau^2;p^m\right)=\tau^2\left(p^m\right),\quad L_1\left(\tau^2;p^m
\right)=p^{22}\tau^2\left(p^{m-1}\right)-2p^{11}\tau(p)\tau\left(p^{m-1}\right)
\tau\left(p^m\right),\nonumber\\
&&L_r\left(\tau^2;p^m\right)=2\left(-p^{11}\right)^r\tau(p)\tau\left(p^{m-1}
\right)\tau\left(p^m\right),\quad r\geq 2\;.\label{r540}
\eea
\end{proposition}
{\sf Proof} $\;\;\;$Prove Proposition by induction. First, calculating the five 
first expressions of $L_r\left(\tau^2;p^m\right)$ one can verify that 
(\ref{r540}) holds. Let Proposition holds for $2\leq q\leq r$, then prove that 
$$
L_{r+1}\left(\tau^2;p^m\right)=2\left(-p^{11}\right)^{r+1}\tau(p)\tau\left(
p^{m-1}\right)\tau\left(p^m\right)\;.
$$
Keeping in mind (\ref{b5}) calculate $L_{r+1}\left(\tau^2;p^m\right)$,
\bea
L_{r+1}\left(\tau^2;p^m\right)\!=\!\tau^2\left(p^{m+r+1}\right)\!-\!\tau^2\left(
p^{r+1}\right)\tau^2\left(p^m\right)\!-\!\tau^2\left(p^r\right)p^{11}\tau\left(
p^{m-1}\right)\left[p^{11}\tau\left(p^{m-1}\right)\!-\!2\tau(p)\tau\left(p^m
\right)\right]\nonumber\\
-2\tau(p)\tau\left(p^{m-1}\right)\tau\left(p^m\right)p^{22}\left[\tau^2\left(
p^{r-1}\right)-\tau^2\left(p^{r-2}\right)p^{11}+\ldots+(-1)^{r+1}\tau^2(p)
p^{11(r-2)}\right]\;.\nonumber
\eea
By (\ref{r532}) the four first terms in the above equality are reduced up to a 
single term,
\bea
\tau^2\left(p^{m+r+1}\right)-\tau^2\left(p^{r+1}\right)\tau^2\left(p^m\right)-
\tau^2\left(p^r\right)p^{11}\tau\left(p^{m-1}\right)\left[p^{11}\tau\left(p^{
m-1}\right)-2\tau(p)\tau\left(p^m\right)\right]=\nonumber\\
2\tau\left(p^m\right)\tau\left(p^{m-1}\right)\tau\left(p^r\right)p^{11}
\left[\tau\left(p^r\right)\tau(p)-\tau\left(p^{r+1}\right)\right]=2\tau\left(
p^m\right)\tau\left(p^{m-1}\right)\tau\left(p^r\right)\tau\left(p^{r-1}\right)
p^{22}\;,\nonumber
\eea
which simplifies further calculations,
\bea
L_{r+1}\left(\tau^2;p^m\right)=2\tau\left(p^m\right)\tau\left(p^{m-1}\right)
p^{22}\;A_{r-1}\;,\quad\mbox{where}\hspace{3cm}\label{r541}\\
A_{r-1}=\tau\left(p^{r-1}\right)\left[\tau\left(p^r\right)-\tau(p)\tau\left(
p^{r-1}\right)\right]+\tau(p)\tau^2\left(p^{r-2}\right)p^{11}-\ldots\pm
\tau^3(p)p^{11(r-2)}\;.\nonumber
\eea
Performing calculations in curl brackets by (\ref{r532}), we obtain 
\bea
L_{r+1}\left(\tau^2;p^m\right)=-2\tau\left(p^m\right)\tau\left(p^{m-1}\right)
p^{33}\;A_{r-2}\;,\quad\mbox{where}\hspace{3cm}\label{r542}\\
A_{r-2}=\tau\left(p^{r-2}\right)\left[\tau\left(p^{r-1}\right)-\tau(p)\tau
\left(p^{r-2}\right)\right]+\tau(p)\tau^2\left(p^{r-3}\right)p^{11}-\ldots\pm
\tau^3(p)p^{11(r-3)}\;.\nonumber
\eea
By comparison (\ref{r541}) and (\ref{r542}) and continuing to contract the terms
in curl brackets, we get
\bea
L_{r+1}\left(\tau^2;p^m\right)=(-1)^r2\tau\left(p^m\right)\tau\left(p^{m-1}
\right)p^{11r}\;A_1\;,\quad A_1=\tau(p)\left[\tau\left(p^2\right)-\tau^2(p)
\right]=-\tau(p)\;p^{11}\;.\nonumber
\eea
Thus, Proposition is proven.$\;\;\;\;\;\;\Box$

Calculate the renormalization function $R_{\infty}\left(\tau^2,p^m\right)$. In 
accordance with Table 1 we obtain $\tau^2\in{\mathbb B}\left\{N^{12}\right\}$, 
i.e., $a_1=12$. Then by Theorem \ref{thm1} and Proposition \ref{pro3} we get
\bea
R_{\infty}\left(\tau^2,p^m\right)=\sum_{r=0}^{\infty}\frac{L_r\left(f;p^m\right)
}{p^{12r}}=\tau^2\left(p^m\right)+p^{10}\tau^2\left(p^{m-1}\right)-\frac{2\;\tau
(p)}{p+1}\tau\left(p^{m-1}\right)\tau\left(p^m\right)\;,\label{r543}
\eea
such that $R_{\infty}\left(\tau^2,1\right)=1$. In the case $m=1$, we find 
$R_{\infty}\left(\tau^2,p\right)=p^{10}+\tau^2(p)\cdot (p-1)/(p+1)$.

Finally, calculate the renormalization function $R_{\infty}\left(\tau^2\cdot 
k^{-25/2},p^m\right)$ such that in accordance with Table 2 we have $\tau^2\cdot 
k^{-25/2}\in{\mathbb B}\left\{N^0\right\}$. Making use of relationship 
(\ref{n18b}) between characteristic functions $L_r\left(f\cdot k^{-s};p^m
\right)$ and $L_r\left(f;p^m\right)$ build the renormalization function as it 
was done in formula (\ref{r538a}) for the Ramanujan $\tau$ Dirichlet series,
\bea
R_{\infty}\left(\frac{\tau^2}{k^{25/2}},p^m\right)=\sum_{r=0}^{\infty}L_r\left(
\frac{\tau^2}{k^{25/2}};p^m\right)=p^{-25m/2}\sum_{r=0}^{\infty}\frac{L_r\left(
\tau^2;p^m\right)}{p^{25r/2}}\;.\nonumber
\eea
By Proposition \ref{pro3} and the above formula we get
\bea
p^{25m/2}R_{\infty}\left(\frac{\tau^2}{k^{25/2}},p^m\right)=\tau^2\left(p^m
\right)+p^{19/2}\tau^2\left(p^{m-1}\right)-\frac{2\;\tau(p)}{p^{3/2}+1}\tau
\left(p^{m-1}\right)\tau\left(p^m\right)\;.\label{r544}
\eea
\section{Numerical Verification}\label{r6}
\begin{figure}[h!]\centering\begin{tabular}{cc}
\includegraphics[height=5cm,width=7.5cm]{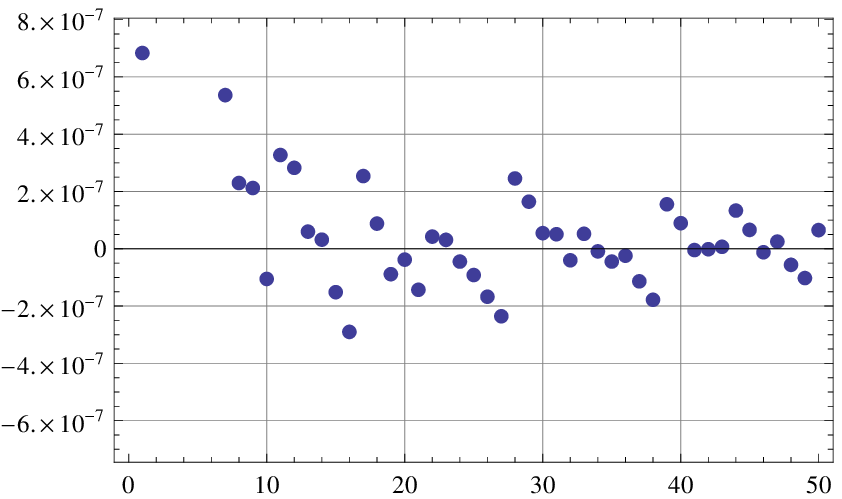} &
\includegraphics[height=5cm,width=7.5cm]{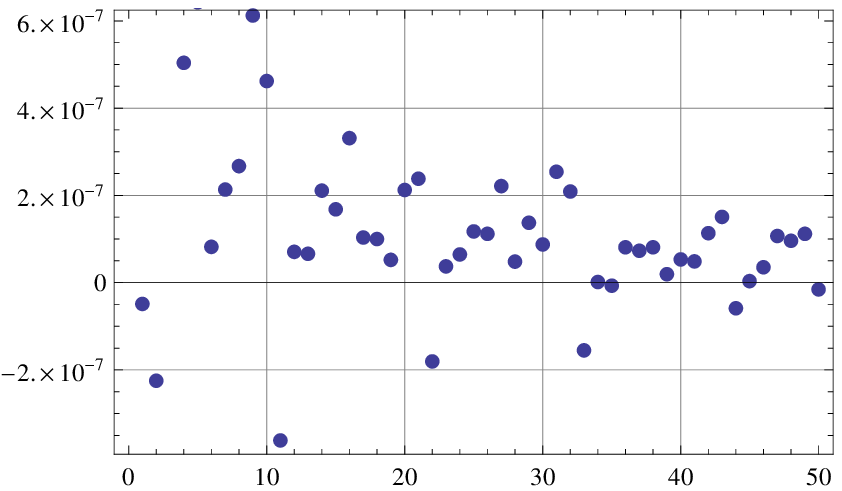}\\
(a) & (b)\end{tabular}
\caption{Plots of $\rho(\phi;N,11)$ (a) and $\rho(\sigma_1;N,11)$ (b) inthe
range $2\cdot 10^4\leq N\leq 10^6$.}\label{fig1}
\end{figure}

\begin{figure}[h!]\centering\begin{tabular}{cc}
\includegraphics[height=5cm,width=7.5cm]{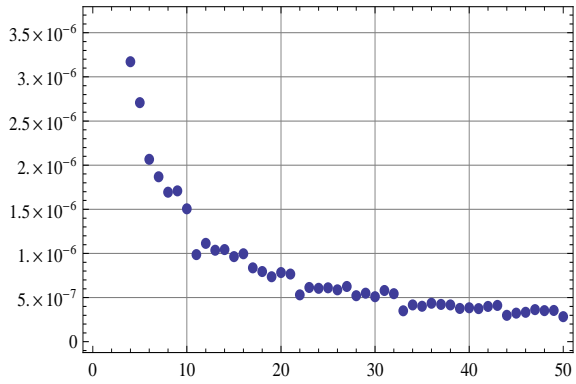} &
\includegraphics[height=5cm,width=7.5cm]{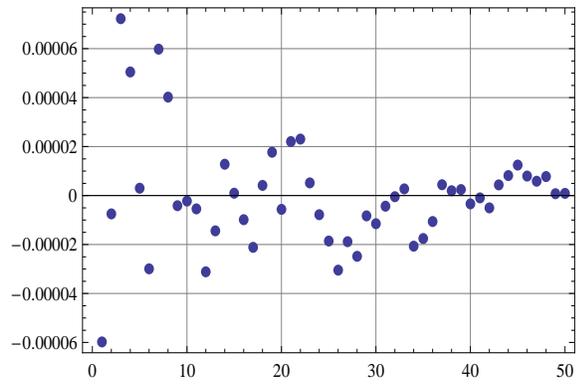}\\
(a) & (b)\end{tabular}
\caption{Plots of $\rho(\sigma_{-1};N,11)$ (a) and $\rho(|\mu|;N,11)$ (b)in
the range $2\cdot 10^4\leq N\leq 10^6$.}\label{fig2}
\end{figure}

\begin{figure}[h!]\centering\begin{tabular}{cc} 
\includegraphics[height=5cm,width=7.5cm]{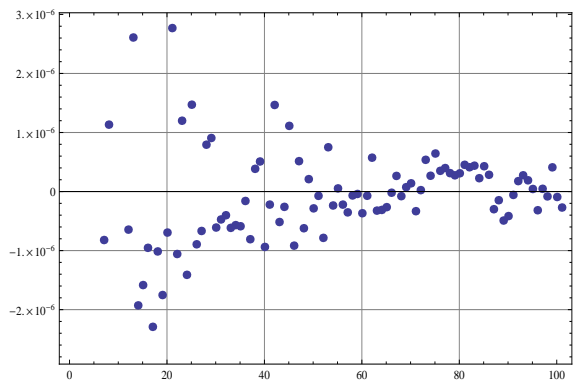} &
\includegraphics[height=5cm,width=7.5cm]{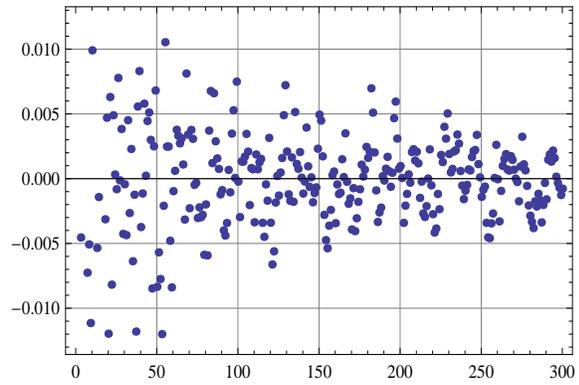}\\
(a) & (b)\end{tabular}
\caption{Plots of $\rho(\tau\cdot k^{-7};N,3)$ in the range $5\cdot 10^3\leq N
\leq 10^5$ (a) and $\rho(\tau^2;N,3)$ (b) in the range $5\cdot 10^3\leq N\leq
3\cdot 10^5$.}\label{fig3}
\end{figure}

In this section we verify the renormalization approach for summation of 
multiplicative arithmetic functions with scaled summation variable by numerical 
calculations. Consider a relative deviation $\rho\left(f;N,p^m\right)$ between 
rescaled $F\left\{f;N,p^m\right\}$ and nonscaled $F\left\{f;N,1\right\}$ 
summatory functions,
\begin{eqnarray}
\rho\left(f;N,p^m\right)=R_{\infty}^{-1}\left(f;p^m\right)\cdot\frac{
F\left\{f;N,p^m\right\}}{F\{f;N,1\}}-1\;,\label{r61}
\end{eqnarray}
where the renormalization function $R_{\infty}\left(f;p^m\right)$ is calculated 
according to Theorem \ref{thm1}. In Figures \ref{fig1}, \ref{fig2} and 
\ref{fig3} we present plots of $\rho\left(f;N,p^m\right)$ for six different 
arithmetic functions. These Figures show that formulas for renormalization 
functions work with high precision.
\section*{Acknowledgement}\label{a8}
The useful discussions with A. Juhasz, A. Mann and Z. Rudnick are highly 
appreciated.

\end{document}